\documentclass[11pt]{article}
\usepackage{graphicx}
\usepackage{color}
\usepackage{amsmath}
\usepackage{amssymb}
\usepackage{amscd}
\usepackage{framed}
\usepackage{bbm}

\usepackage{fancyhdr,a4wide}
\usepackage{amsthm}

\newcommand{\R}{\mathbb{R}}
\newcommand{\inr}[1]{\left\langle #1 \right\rangle}
\newcommand{\sgn}{{\rm sgn}}

\newcommand{\E}{\mathbb{E}}

\newcommand{\eps}{\varepsilon}

\newtheorem{Theorem}{Theorem}[section]
\newtheorem{Lemma}[Theorem]{Lemma}

\newtheorem{Definition}[Theorem]{Definition}

\newtheorem{Corollary}[Theorem]{Corollary}
\newtheorem{Remark}[Theorem]{Remark}
\newtheorem{Example}[Theorem]{Example}
\newtheorem{Assumption}{Assumption}[section]

\numberwithin{equation}{section}

\def \proof {\noindent {\bf Proof.}\ \ }

\def \endproof
{{\mbox{}\nolinebreak\hfill\rule{2mm}{2mm}\par\medbreak}}
\def\IND{\mathbbm{1}}

\def\IND{\mathbbm{1}}

\begin{document}
\title{An isomorphic Dvoretzky-Milman Theorem using general random ensembles}
\author{
Shahar Mendelson \thanks{Mathematical Sciences Institute, The Australian National University. Email: shahar.mendelson@anu.edu.au} }

\maketitle

\begin{abstract}
We construct rather general random ensembles that yield the optimal (isomorphic) estimate in the Dvoretzky-Milman Theorem. This is the first construction of non gaussian/spherical ensembles that exhibit the optimal behaviour. The ensembles constructed here need not satisfy any rotation invariance and can be rather heavy-tailed.
\end{abstract}

\section{Introduction}
The Dvoretzky-Milman Theorem is one of the cornerstones of Asymptotic Geometric Analysis. The original result, due to Dvoretzky (\cite{MR0139079}) is that $\ell_2$ is \emph{finitely represented} in any infinite dimensional normed space. That means that given an infinite dimensional normed space $(E,\| \ \|)$, for any $\eps>0$ and integer $k$ there is a subspace $E_k \subset E$ of dimension $k$ such that $d(\ell_2^k,E_k) \leq 1+\eps$. Here $\ell_2^k=(\R^k, \| \ \|_2)$ and the distance is the Banach-Mazur distance, defined by
$$
d(E,F) = \inf \left\{ \|T\| \|T^{-1}\| \ : \ T:E \to F \ {\rm is \ invertible} \ \right\}.
$$
Dvoretzky's result was actually quantitative, but with a suboptimal estimate. Milman, in his seminal work \cite{MR0293374} (see also \cite{MR3331351} for an extensive exposition), used the idea of concentration of measure, and specifically, the concentration of Lipschitz functions on the Euclidean sphere around their medians, to obtain the optimal finite dimensional estimate. Milman showed that each convex body\footnote{A convex body in $\R^n$ is a convex, centrally symmetric set with a nonempty interior. There is an obvious equivalence between a convex body $K$ and the norm $\| \ \|_K$ whose unit ball is $K$.} $K \subset \R^n$ has a critical dimension $d^*(K)$ and for every $\eps>0$, a typical (with respect to the Haar measure on the grassmannian) subspace of $\R^n$ of dimension $s=c(\eps)d^*(K)$ satisfies that $d(\ell_2^s, K \cap E) \leq 1+\eps$. Moreover, he showed  that there is a linear image $K^\prime$ of $K$ for which $d^*(K^\prime) \geq c \log n$ for an absolute constant $c$.

The following is the equivalent gaussian formulation of Milman's result, due to Pisier (\cite{MR864714}, see also \cite{MR1036275}).

Let $K \subset \R^n$ be a convex body and denote by $\| \ \|_K$ the norm on $\R^n$ whose unit ball is $K$. Let $K^\circ$ be the dual body, i.e., for every $x \in \R^n$, $\|x\|_K = \sup_{t \in K^\circ} \inr{x,t}$. Let $(g_i)_{i=1}^n$ be independent standard gaussian random variables and set
$$
\ell(K)=\E \left\|\sum_{i=1}^n g_i e_i \right\|_K = \E \sup_{t \in K^\circ} \sum_{i=1}^n g_i t_i.
$$
\begin{Definition} \label{def:critical}
 The critical dimension of $K$ is
$$
d^*(K)=\left(\frac{\ell(K)}{\sup_{t \in K^\circ} \|t\|_2}\right)^2.
$$
\end{Definition}

Here and in what follows, we denote by $B_2^n$ the unit ball of $\ell_2^n$, and $S^{n-1}$ is the Euclidean sphere in $\R^n$.

The gaussian version of Milman's result is as follows:

\begin{Theorem} \label{thm:DM-intro}
For every $\eps>0$ there exist constants $c_1$ and $c_2$ that depend on $\eps$ such that the following holds. Let $K \subset \R^n$ be a convex body and set $d = c_1(\eps)d^*(K)$. Let $T_g:\R^d \to \R^n$ be the random operator whose entries are independent, standard gaussian random variables and set $E = T_g \R^d$. Then with probability at least $1-2\exp(-c_2(\eps) d^*(K))$,
$$
(1-\eps)\ell(K) (K \cap E) \subset T_g B_2^d \subset (1+\eps) \ell(K) (K \cap E)
$$
\end{Theorem}

\begin{Remark}
For reasons that will become clear immediately, the dependence on $\eps$ is not of importance in the context of this article. However, it is a question that has been studied extensively over the years, and we refer the reader to the book \cite{MR3331351} for a detailed discussion on that topic.
\end{Remark}

Note that Theorem \ref{thm:DM-intro} gives much more than the existence of an almost Euclidean subspace of $(\R^n,\| \ \|_K)$. Rather, it shows that a \emph{typical subspace} (of the right dimension) with respect to the gaussian ensemble is almost Euclidean; Milman's original proof shows that the same phenomenon is true with respect to the spherical ensemble.

To explain the role of the gaussian/spherical ensemble in Theorem \ref{thm:DM-intro} and highlight the difficulty in obtaining non-gaussian/non-spherical versions of that result, let us present an outline of the proof of Theorem \ref{thm:DM-intro}.

\vskip0.3cm

The proof of Theorem \ref{thm:DM-intro} is remarkably simple: an operator $T:\ell_2^d \to (\R^n,\| \ \|_K)$ is an almost isometric embedding of $\ell_2^d$ in $(\R^n, \| \ \|_K)$ if there is some $\alpha>0$ such that for every $x \in S^{d-1}$,
\begin{equation} \label{eq:gaussian-intro}
(1-\eps)\alpha \leq \left\|\sum_{i=1}^d x_i Te_i\right\|_K \leq (1+\eps) \alpha;
\end{equation}
it is an isomorphic embedding if $1 \pm \eps$ is replaced by absolute constants $c$ and $C$. When $T$ is a gaussian operator, each $Te_i$ is distributed as $G$, the standard gaussian vector in $\R^n$, and the $d$ random vectors $\{Te_i, 1 \leq i \leq d\}$ are independent. Thanks to rotation invariance, for each $x \in S^{d-1}$, the random vector $\sum_{i=1}^d x_i Te_i$ has the same distribution as $\|x\|_2 G=G$, and by the gaussian concentration theorem (see, e.g. \cite{MR1036275}), there is an absolute constant $c$ such that for any $u>0$,
\begin{equation} \label{eq:gaussian-conc-intro}
Pr\left( \left| \|G\|_K - \ell(K) \right| \geq u \right) \leq 2\exp\left(-c\frac{u^2}{\sup_{t \in K^\circ} \|t\|_2^2}\right).
\end{equation}
Setting $u=\eps \ell(K)$, it follows from \eqref{eq:gaussian-conc-intro} that for any $x \in S^{d-1}$, \eqref{eq:gaussian-intro} holds with probability at least $1-2\exp(-c\eps^2d^*(K))$ for $\alpha=\ell(K)$; the high probability estimate allows one to have uniform control on an $\eps$-net in $S^{d-1}$; and thanks to standard convexity and approximation arguments the uniform estimate on $S^{d-1}$ can be ensured. We refer the reader to \cite{MR1036275} for the complete proof.

\vskip0.3cm

While this proof is undeniably beautiful, it is also rather restrictive. It is based on two crucial features of the gaussian vector: that for every $x \in S^{d-1}$ the random variable $\left\|\sum_{i=1}^d x_i Te_i\right\|_K$ concentrates around its mean, and that, by rotation invariance, the means $\E \left\|\sum_{i=1}^d x_i Te_i\right\|_K$ are the same (or, for an isomorphic estimate, are equivalent). These features are shared by the spherical ensemble.

\begin{framed}
In what follows, we refer (somewhat inaccurately) to a random ensemble $T$ as \emph{an optimal Dvoretzky-Milman ensemble} if it exhibits the fact that an arbitrary convex body $K$ has a subspace of dimension proportional to $d^*(K)$ that is isomorphic (with absolute constants) to a Euclidean space. More accurately, an optimal Dvoretzky-Milman ensemble $T$ satisfies that for every convex body $K \subset \R^n$, with probability at least $1/2$, for $d=c_0 d^*(K)$ and every $x \in \R^d$,
\begin{equation} \label{eq:DM-ensemble}
c\ell(K)\|x\|_2 \leq \|Tx\|_K \leq C\ell(K)\|x\|_2.
\end{equation}
Here, $c_0, c$ and $C$ are constants that might depend on properties of the ensemble but are independent of $K$, $n$ or $d^*$.
\end{framed}

\vskip0.3cm
As surprising as it may seem, it is hard to obtain a behaviour like \eqref{eq:DM-ensemble} when leaving the gaussian/spherical realm, even if the ensemble is generated by independent copies of a very well behaved random variable---as the next example shows.

\begin{Example} \label{ex:Bernoulli-bad}
Let $K=B_\infty^n$, the unit cube in $\R^n$; its dual body is $B_1^n$, the unit ball in $\ell_1^n =(\R^n, \| \ \|_1)$. It is standard to verify that
$$
\ell(B_\infty^n)=\E \max_{1 \leq i \leq n} |g_i| \sim \sqrt{\log n},
$$
and clearly, $\sup_{t \in B_1^n} \|t\|_2 =1$. Therefore, the critical dimension of the unit cube is $d^*(B_\infty^n) \sim \log n$. By Theorem \ref{thm:DM-intro}, there are absolute constants $\delta>0$, $c_0,c$ and $C$, such that for $d = c_0 \log n$, with probability at least $1-2/n^\delta$,
$$
c\sqrt{\log n} \leq \inf_{x \in S^{d-1}} \left\|T_gx\right\|_\infty \leq \sup_{x \in S^{d-1}} \left\|T_gx\right\|_\infty \leq C\sqrt{\log n}.
$$
Thus, on that event, the random function $\psi(x)=\|T_g x\|_K$ is equivalent to $\sqrt{\log n}$ on $S^{d-1}$.

On the other hand, let $\zeta$ be uniformly distributed in $[-1,1]$. Set
$$
T_\zeta=\left\{(\zeta_{ij}) : 1 \leq i \leq  n, \ 1 \leq j \leq d\right\}
$$
to be the random operator whose entries are independent copies of $\zeta$.

Note that here we no longer have rotation invariance, and it turns out that with high probability, the random function
$$
\phi(x)=\|T_\zeta x\|_\infty=\left\|\sum_{j=1}^d x_j T_\zeta e_j \right\|_\infty
$$
is not equivalent to one value on $S^{d-1}$. Indeed, for $y_1=(1,0,...0)$, $\phi(y) \leq 1$ almost surely. At the same time, for $\eta=(\eta_j) \in \{-1,1\}^d$ let $y_\eta=(\eta_j/\sqrt{d})_{j=1}^d \in S^{d-1}$, and
$$
\phi(y_\eta) = \frac{1}{\sqrt{d}} \max_{i=1,...n} \left|\sum_{j=1}^d \eta_j \zeta_{ij}\right|.
$$
Denote by $|A|$ the cardinality of a set $A$ and observe that
\begin{equation} \label{eq:prob-est-zeta}
Pr \left(\{ \exists i \  : \ |\{j: |\zeta_{ij}| \geq 1/2\}| \geq d/10\} \right) \geq 1-2\exp(-c_1dn)
\end{equation}
for an absolute constant $c_1$. For every realization of the event from \eqref{eq:prob-est-zeta}, let $i^*$ be the corresponding index and set $\eta_j=\sgn(\zeta_{i^*j})$ for $1 \leq j \leq d$.  Clearly,
$$
\phi(y_\eta) = \frac{1}{\sqrt{d}} \sum_{j=1}^d |\zeta_{i^*j}| \geq c_2 \sqrt{d}
$$
for an absolute constant $c_2$, and recalling that $d=c_0\log n$, we have that with probability at least $1-2\exp(-c_0 n \log n)$
$$
\inf_{x \in S^{d-1}} \left\|\sum_{i=1}^d x_i T_\zeta e_i \right\|_\infty \leq 1 \ \ \  {\rm and}  \ \ \  \sup_{x \in S^{d-1}} \left\|\sum_{i=1}^d x_i T_\zeta e_i \right\|_\infty \geq \sqrt{ c_2 \log n}.
$$
\endproof
\end{Example}

This example illustrates a rather surprising fact: although 50 years have passed since Milman's proof, there are no examples of random ensembles that are optimal Dvoretzky-Milman ensembles---other than the gaussian ensemble or the spherical one. The few known ``non gaussian" versions of an isomorphic Dvoretzky-Milman Theorem hold only for convex bodies $K$ that have special structure, like cotype $2$ (see, e.g. \cite{MR2391154}).

\vskip0.3cm

Here we present a rather general construction that yields optimal Dvoretzky-Milman ensembles. The random operators need not satisfy rotation invariance of any kind, nor do they exhibit a gaussian-like concentration; in fact, they can be rather heavy-tailed.

\vskip0.3cm

Before presenting a more accurate description of the construction, let us present two examples. A third example is presented in Section \ref{sec:example-heavy-tailed}.

\begin{Theorem} \label{thm:bernoulli-good}
There exist absolute constants $c_0,c_1,c_2,c$ and $C$ such that the following holds: Let $K \subset \R^n$ and set $d^*(K)$ to be its critical dimension. Let $d=c_0d^*(K)$ and $m = c_1n$, and set
$$
\Gamma_1=m^{-1/2}(\zeta_{ij}):\R^n \to \R^m, \ \ \ \ \Gamma_2=(\zeta^\prime_{ij}):\R^m \to \R^d,
$$
where $\{ (\zeta_{ij}) : 1 \leq i \leq m, \ 1 \leq j \leq n\}$ and $\{ (\zeta_{ij}^\prime) : 1 \leq i \leq d, \ 1 \leq j \leq m\}$ are independent random variables that are uniformly distributed in $[-1,1]$. Define
$$
\Gamma:\Gamma_1^*\Gamma_2^* : \R^d \to \R^n.
$$
Then, with probability at least $1-2\exp(-c_2 d^*(K))$, for every $x \in S^{d-1}$,
$$
c \ell(K) \leq \|\Gamma x\|_K \leq C \ell(K).
$$
\end{Theorem}
Theorem \ref{thm:bernoulli-good} recovers the isomorphic Dvoretzky-Milman Theorem using the product of two matrices of the type $T_\zeta$ of appropriate dimensions. However, as Example \ref{ex:Bernoulli-bad} shows, a \emph{single matrix} of the same type is not a suitable choice.
\vskip0.3cm

Before we formulate the second example, let us recall some standard definitions.
\vskip0.3cm
\begin{Definition} \label{def:subgaussian}
Let $Y$ be a centred random vector in $\R^\ell$.
\begin{description}

\item{$\bullet$} $Y$ is \emph{isotropic} if its covariance is the identity; that is, if for every $t \in \R^\ell$
$$
\E \inr{Y,t}^2 = \|t\|_2^2.
$$
\item{$\bullet$} $Y$ is \emph{log-concave} if has a density that is a log-concave function.
\item{$\bullet$} $Y$ is $L$-subgaussian if for any $t \in \R^\ell$ and every $p \geq 2$,
$$
\|\inr{Y,t}\|_{L_p} \leq L \sqrt{p} \|\inr{Y,t}\|_{L_2}.
$$
\item{$\bullet$} Let $q > 2$. $Y$ satisfies an $L_q-L_2$ norm equivalence with constant $L$ if for every $t \in \R^\ell$,
$$
\|\inr{Y,t}\|_{L_q} \leq L \|\inr{Y,t}\|_{L_2}.
$$
\end{description}
\end{Definition}
Observe that if, in addition to being $L$-subgaussian, $Y$ is isotropic, then for every $t \in \R^\ell$ and $p \geq 1$, $\|\inr{Y,t}\|_{L_p} \leq L \sqrt{p}\|t\|_2$. Clearly, the analogous observation is true if $Y$ satisfies $L_q-L_2$ norm equivalence with constant $L$.
\begin{Remark}
To put these definitions in some context, the random vector $Y=(\zeta_i)_{i=1}^\ell$ is isotropic; $L$-subgaussian for an absolute constant $L$; and log-concave.
\end{Remark}

\begin{framed}
In what follows, all the random vectors we consider are symmetric and isotropic. Given integers $n$ and $d$, we consider random vectors $Z$ in $\R^n$ and random vectors $X$ in $\R^d$. For $m \geq n$ let  $Z_1,...,Z_m$ be independent copies of $Z$ and set $X_1,...,X_m$ to be independent copies of $X$. Define $\Gamma_1$ to be the random matrix whose rows are $Z_i/\sqrt{m}, 1 \leq i \leq m$, and set $\Gamma_2$ to be the matrix whose columns are $X_1,...X_m$. Finally, let
\begin{equation} \label{eq:Gamma-intro}
\Gamma=\Gamma_1^* \Gamma_2^*:\R^d \to \R^n.
\end{equation}
As we explain, if the rows of $\Gamma_1$ and the columns of $\Gamma_2$ satisfy an appropriate mixture of tail and small-ball properties then $\Gamma=\Gamma_1^*\Gamma_2^*$ is an optimal Dvoretzky-Milman ensemble.
\end{framed}
The second outcome of our main result is as follows:
\begin{Theorem} \label{thm:log-concave}
For every $L$, $L^\prime \geq 1$ there are constants $c_0,...,c_2$, $c$ and $C$ that depend on $L$ and $L^\prime$ such that the following holds. Consider a convex body $K \subset \R^n$, set $d=c_0 d^*(K)$ and put $m=c_1n$.
Let $Z$ be an isotropic, $L$-subgaussian random vector in $\R^n$ and let $X$ be an isotropic, log-concave random vector in $\R^d$. Assume that the densities of all of the one-dimensional marginals $\inr{Z,t}$ and $\inr{X,v}$ are bounded by $L^\prime$.

If $\Gamma$ is as in \eqref{eq:Gamma-intro}, then with probability at least
$$
1-2\exp\left(-c_2\min\{(n d^*(K))^{1/4},d^*(K)\}\right),
$$
for every $x \in S^{d-1}$,
$$
c\ell(K) \leq \left\| \Gamma x \right\|_K \leq C \ell(K).
$$
\end{Theorem}
\begin{Remark}
Although $Z$ and $X$ need not have iid coordinates, they still share some of the features of $(\zeta_i)$. In particular, the one-dimensional marginals of $X$ exhibit a fast tail decay. We shall present another example in Section \ref{sec:example-heavy-tailed} which shows that $X$ can be heavy-tailed.
\end{Remark}

Let us turn to the formulation of the main result of this note, beginning with the required features of $Z$ and $X$.

\subsection{The Assumptions}
The first assumption we require is that both $X$ and $Z$ satisfy a rather weak small-ball property:
\begin{Assumption} \label{ass:small-ball-assumption}
There is a constant $\kappa_0$ such that for every $t \in \R^n$ and every $v \in \R^d$,
$$
Pr\left( |\inr{Z,t}| \leq \kappa_0 \|t\|_2\right) \leq \frac{1}{1000} \ \ \ {\rm and} \ \ \ Pr\left( |\inr{X,v}| \leq \kappa_0 \|v\|_2\right) \leq \frac{1}{1000}.
$$
\end{Assumption}
\begin{Remark}
The choice of $1/1000$ is completely arbitrary; in fact, any uniform probability estimate strictly smaller than $1/2$ suffices for our purposes.
\end{Remark}

The second assumption deals with the tail behaviour of one dimensional marginals of $X$ and $Z$.

\begin{Assumption} \label{ass:on-X}
Assume that
\begin{description}
\item{$(1)$} $Z$ is $L$-subgaussian.
\item{$(2)$} $X$ satisfies an $L_q-L_2$ norm equivalence for some $q>2$ with constant $L$. In particular, since $X$ is isotropic, it follows that for every $v \in \R^d$, $\|\inr{X,v}\|_{L_q} \leq L\|v\|_2$.
\end{description}
\end{Assumption}

In itself, Assumption \ref{ass:on-X} is not enough and additional information on $\Gamma_2$ is required. To formulate this final assumption, denote by $\|a\|_0$ the cardinality of the support of the vector $a$, i.e., $\|a\|_0=|\{i: a_i \not= 0\}|$.
\begin{Assumption} \label{ass:on-Gamma-2}
For $\kappa_1 \geq 1$ and well chosen $0<\delta,\theta<1/4$, there is a nontrivial event ${\cal A}$ such that for any $(X_i)_{i=1}^m \in {\cal A}$,
\begin{equation} \label{eq:operator-norm}
\sup_{v \in B_2^d} \|\Gamma_2 v\|_2 \leq \kappa_1 \sqrt{m},
\end{equation}
and
\begin{equation} \label{eq:rearrangement}
\sup_{\{a \in \R^m : \|a\|_2 \leq 1, \ \|a\|_{0} \leq \theta m\}} \|\sum_{i=1}^m a_i X_i\|_2 \leq \delta \sqrt{m}.
\end{equation}
\end{Assumption}

\begin{Remark}
Assumption \ref{ass:on-Gamma-2} has strong ties to the study of the extremal singular values of $\Gamma_2$. Thanks to these connections, there are many scenarios in which \eqref{eq:operator-norm} and \eqref{eq:rearrangement} are known to be true on a large event. Indeed, \eqref{eq:operator-norm} implies that the largest singular value of $\Gamma_2/\sqrt{m}$ satisfies $\lambda_{\max}(\Gamma_2/\sqrt{m}) \leq \kappa_1$. And in fact, in all the examples we present here, $\kappa_1=2$ suffices: under rather mild assumptions on $X$, for any $m \geq d$, with high probability
$$
\lambda_{\max}(\Gamma_2/\sqrt{m}) \leq 1+c^\prime \sqrt{d/m}.
$$
The second condition in the definition of ${\cal A}$, Equation \eqref{eq:rearrangement}, is closely related to Bai-Yin type estimates on $\lambda_{\max}(\Gamma_2/\sqrt{m})$ and $\lambda_{\min}(\Gamma_2/\sqrt{m})$ (see, e.g., \cite{MR2601042,MR3191978,MR3872323}).
\end{Remark}

Although the event ${\cal A}$ depends only on $(X_i)_{i=1}^m$, we view it with respect to the product $(X_i,Z_i)_{i=1}^m$. Note that the event actually depends on five parameters: $m$, $d$ (because $X$ is a random vector in $\R^d$), $\kappa_1$, $\delta$ and $\theta$. While $\kappa_1$ can be arbitrary (and, as noted previously, in most interesting cases, $\kappa_1=2$ suffices), a suitable choice of $\delta$ and $\theta$ requires more care: for the embedding result to hold, these parameters need to satisfy several constraints that depend on the values of $\kappa_0$, $\kappa_1$, $L$ and $q$. As a result, the event ${\cal A}$ should be understood  as corresponding to a fixed choice of $m$, $d$, $\kappa_1$, $\delta$ and $\theta$; the proof will dictate the constraints on $m$, $d$, $\delta$ and $\theta$, and one has to show that for a set of parameters that satisfy those constraints the event ${\cal A}$ has nontrivial probability. To be more accurate:
\begin{framed}
\begin{description}
\item{$\bullet$}  The proof of our main result is based on a net argument. From here on, fix $0<\rho<1/4$ and denote by $V$ a maximal $\rho$-separated subset of $S^{d-1}$ with respect to $\ell_2$ norm. In what follows, we shall call such a set \emph{a $\rho$-net}.  By a volumetric estimate
$$
|V| \leq \exp(d\log(5/\rho)),
$$
a fact we use frequently.

The right value of $\rho$ turns out to be a constant that depends on $\kappa_0$, $\kappa_1$, $L$ and $q$ (see the proof of Theorem \ref{thm:main-intro}, below).

\item{$\bullet$} Given a suitable choice of $\rho$, $m$ has to satisfy that
\begin{equation} \label{eq:intro-m}
m \geq c_0\max\left\{\frac{d^*(K)}{\rho}, n\right\}
\end{equation}
for a constant $c_0$ that depends on $\kappa_0$, $\kappa_1$, $L$ and $q$.

\item{$\bullet$} The parameters $\delta$ and $\theta$ must satisfy that for suitable constants $c_1$ and $c_2$ that depend on $\kappa_0,\kappa_1,L$ and $q$,
\begin{equation} \label{eq:intro-theta-delta}
\delta \leq \frac{c_1}{\log(5/\rho)}, \ \ \ {\rm and} \ \ \ \theta^{(q-2)/2(q+2)} \sqrt{\log(e/\theta)} \leq \frac{c_2}{\log(5/\rho)}.
\end{equation}
\item{$\bullet$} At the same time, for values of $\delta$ and $\theta$ that satisfy \eqref{eq:intro-theta-delta}, the dimension $d$ of the subspace that is close to Euclidean will be
\begin{equation} \label{eq:intro-d}
d = c_3(\kappa_0,L) (\theta^{4/(2+q)}/\log^5(5/\rho)) d^*(K).
\end{equation}
\end{description}
\end{framed}

\subsection{The main result}
We show that
\begin{Theorem} \label{thm:main-intro}
There are constants $c$, $c^\prime$ and $C$ that depend on $\kappa_0$, $\kappa_1$, $L$ and $q$ such that the following holds. Set $m$ as in \eqref{eq:intro-m}, $\theta$ and $\delta$ as in \eqref{eq:intro-theta-delta} and $d$ as in \eqref{eq:intro-d}. Then with probability at least $Pr({\cal A}) - 2\exp(-c^\prime d^*(K))$, for every $x \in S^{d-1}$
$$
c \ell(K)  \leq \|\Gamma x\| \leq C \ell(K).
$$
\end{Theorem}

\vskip0.3cm

Writing $\inr{\Gamma v ,t}=\inr{\Gamma_2^*v,\Gamma_1 t}$, Theorem \ref{thm:main-intro} follows once we show that there is an event ${\cal B}$ of probability at most $2\exp(-c^\prime d^*(K))$ such that on ${\cal A}\backslash {\cal B}$,
$$
\sup_{v \in B_2^d} \sup_{t \in K^\circ} \inr{ \Gamma_2^* v ,\Gamma_1 t} \leq C\ell(K),
$$
and
$$
\inf_{v \in S^{d-1}} \sup_{t \in K^\circ} \inr{ \Gamma_2^* v ,\Gamma_1 t} \geq c\ell(K).
$$
To that end, let $0<\rho<1/4$ be specified later. By convexity (for the standard argument, see, e.g., \cite{MR3331351}), it is enough to establish the upper estimate on $U$ that is a $1/5$-net in $B_2^d$; and by a volumetric estimate, there is a $1/5$-net in $B_2^d$ whose cardinality is at most $\exp(10d)$.

Thus, for the upper bound it suffices to show that:
\begin{Theorem} \label{thm:upper}
There are constants $c$, $c^\prime$ $C_0$ and $C$ that depend only on $\kappa_1$, $L$ and $q$ for which the following holds. Let $d \leq n$ and $m \geq C_0n$. There is an event ${\cal B}_1$ of probability at most  $2\exp(-c^\prime d)$ such that on ${\cal A} \backslash {\cal B}_1$, for every $v \in U$,
$$
\sup_{t \in K^\circ} \inr{ \Gamma_2^* v ,\Gamma_1 t} \leq C\ell(K).
$$
\end{Theorem}

The reverse inequality is far more delicate. The key step is to obtain a uniform lower bound on the $\rho$-net $V \subset S^{d-1}$.
\begin{Theorem} \label{thm:lower}
There are constants $c_0,...,c_4$ that depend on $\kappa_0,\kappa_1$, $L$ and $q$ such that the following holds. Fix $0<\delta,\theta<1/4$, let
$$
m \geq c_0\max\{d^*(K)/\rho, n\} \ \ \ {\rm and} \ \ \ d = c_1 (\theta^{4/(2+q)}/\log^5(5/\rho)) d^*(K).
$$

There is an event ${\cal B}_2$ such that $Pr({\cal B}_2) \leq 2\exp(-c_2 d \log(5/\rho))$ and on ${\cal A}\backslash {\cal B}_2$, for every $v \in V$,
$$
\sup_{t \in K^\circ} \inr{ \Gamma_2^* v ,\Gamma_1 t} \geq \left(c_3-c_4 \log(5/\rho) \left(\delta+\theta^{(q-2)/2(q+2)}\sqrt{\log(e/\theta)}\right) \right) \ell(K) .
$$
In particular, if
$$
\delta \leq \frac{c_3}{8c_4\log(5/\rho)}, \ \ \ {\rm and} \ \ \ \theta^{(q-2)/2(q+2)} \sqrt{\log(e/\theta)} \leq \frac{c_3}{8c_4\log(5/\rho)},
$$
then on the same event, for every $v \in V$,
$$
\sup_{t \in K^\circ} \inr{ \Gamma_2^* v ,\Gamma_1 t} \geq \frac{c_3}{2}\ell(K).
$$
\end{Theorem}

\vskip0.3cm
The combination of Theorem \ref{thm:upper} and Theorem \ref{thm:lower} leads to the proof of Theorem \ref{thm:main-intro}:

\vskip0.3cm
\noindent{\bf Proof of Theorem \ref{thm:main-intro}.} Since the conditions on $d$ and $m$ imposed in Theorem \ref{thm:lower} are more restrictive than those in Theorem \ref{thm:upper}, we choose $d$ and $m$ as in the latter. All that remains is to prove the uniform lower bound and determine a suitable choice of $\rho$.

Let ${\cal B}_1$ and ${\cal B}_2$ be the events from Theorem \ref{thm:upper} and Theorem \ref{thm:lower}, respectively.

For every $x \in S^{d-1}$ let $v \in V$ satisfy that $\|x-v\|_2 \leq \rho$. On ${\cal A} \backslash ({\cal B}_1 \cup {\cal B}_2)$, we have, using the same notation as in the two theorems,
\begin{align*}
\sup_{t \in K^\circ} \inr{\Gamma_1^* \Gamma_2^* x ,t} \geq & \sup_{t \in K^\circ} \inr{ \Gamma_2^* v ,\Gamma_1 t} - \sup_{t \in K^\circ} \inr{ \Gamma_2^* (x-v) ,\Gamma_1 t}
\\
\geq & \left(\frac{c_3}{2}-C \rho\right) \ell(K) \geq \frac{c_3}{4} \ell(K),
\end{align*}
provided that $\rho \leq c_3/4C$. Since $c_3$ and $C$ depend on $\kappa_0$, $\kappa_1$, $L$ and $q$, so does $\rho$; and as a result so do $\theta$, $\delta$, $m$ and $d$ with the choices specified in \eqref{eq:intro-theta-delta}, \eqref{eq:intro-m} and \eqref{eq:intro-d}.
\endproof

In Section \ref{sec:upper} we prove Theorem \ref{thm:upper} and Theorem \ref{thm:lower} is proved in Section \ref{sec:lower}. The proofs are based on Talagrand's generic chaining mechanism and properties of Bernoulli processes that are outlined in Section \ref{sec:chaining}. A key component of the proof of Theorem \ref{thm:lower} is that the suprema of certain Bernoulli processes, indexed by random subsets of $\R^m$, actually dominate the expected suprema of their gaussian counterparts, and do so with high probability. This type of domination goes in the \emph{opposite direction} of what one usually expects. It holds here only because of the special structure of the (random) indexing sets, and the role of $\Gamma_1$ and $\Gamma_2$ is to generate that special structure with sufficiently high probability.

Section \ref{sec:examples} is devoted to the proofs of Theorem \ref{thm:bernoulli-good}, Theorem \ref{thm:log-concave} and to another example, presented in Section \ref{sec:example-heavy-tailed}. The latter shows that $X$ can be heavy-tailed, but still $\Gamma$ is an optimal Dvoretzky-Milman ensemble.

\vskip0.3cm

We end this introduction with a word about notation. Throughout, absolute constants are denoted by $c$,$C$,$c^\prime$, etc. Their value may change from line to line, though they remain unchanged within each proof/statement. Unless mentioned otherwise, all constants are absolute---meaning that they are just positive numbers that do not depend on any of the parameters of the problem. If a constant does depend on some parameter, that will be made explicit by writing $c(L)$, $c(L,q)$, etc. We write $a \lesssim b$ if there is an absolute constant $c$ such that $a \leq cb$, and $a \sim b$ means that $a \lesssim b$ and $b \lesssim a$.

\section{Examples} \label{sec:examples}
Let us return to the two examples mentioned previously. The key is to show that in both cases ${\cal A}$ is a large event.
\subsection{Proof of Theorem \ref{thm:bernoulli-good}.}
Recall that $\zeta$ is distributed uniformly in $[-1,1]$ and that $X=(\zeta_i)_{i=1}^d$ and $Z=(\zeta^\prime_i)_{i=1}^n$ have independent coordinates distributed as $\zeta$. Observe that $X$ and $Z$ are isotropic random vectors that are $L$-subgaussian for an absolute constant $c$. Also, it is standard to verify that they satisfy the small ball property of Assumption \ref{ass:small-ball-assumption}. Therefore, both $\kappa_0$ and $L$ are absolute constants in this case.

All that remains is to show that ${\cal A}$ has high probability when $X$ is an arbitrary $L$-subgaussian random vector. As a result, the assertion of Theorem \ref{thm:bernoulli-good} actually holds when $X$ and $Z$ are subgaussian random vectors that satisfy the small-ball property. In such a case, for $c=c(\kappa_0,L)$, with probability at least $1-2\exp(-cd^*(K))$, $\Gamma$ is an optimal Dvoretzky-Milman ensemble.

Before turning to the proof, let us describe certain features of subgaussian random variables. The proofs of all these facts are standard and can be found, for example, in \cite{MR1102015,MR1385671,MR3113826}.
\begin{Definition} \label{def-psi-alpha}
Let $\xi$ be a centred random variable. For $1 \leq \alpha \leq 2$, the $\psi_\alpha$ norm of $\xi$ is defined by
$$
\inf \left\{ c>0 : \E \exp(|\xi/c|^\alpha) \leq 2\right\}.
$$
\end{Definition}

There is a well-known equivalence between the $\psi_\alpha$ norm of a random variable, the growth of its moment and its tail behaviour:
\begin{Theorem} \label{thm:psi}
Let $1 \leq \alpha \leq 2$ and consider a centred random variable $\xi$. The following are equivalent:
\begin{description}
\item{$(1)$} $\|\xi\|_{\psi_\alpha} \leq L$;
\item{$(2)$} For every $p \geq 2$, $\|\xi\|_{L_p} \leq L^\prime p^{1/\alpha}$;
\item{$(3)$} For every $u \geq 1$, $Pr( |\xi| \geq L^{\prime \prime} u) \leq 2\exp(-u^\alpha)$.
\end{description}
Moreover, the constants $L$, $L^\prime$ and $L^{\prime \prime}$ are equivalent in the sense that $L \sim L^\prime \sim L^{\prime \prime}$.
\end{Theorem}
Next, let us turn to a useful feature of independent random variables with bounded $\psi_2$ norms:
\begin{Lemma} \label{lemma-psi-2-tensor}
There exists an absolute constant $c$ such that the following holds. Let $\xi_1,...,\xi_\ell$ be independent, centred random variables. Then for any $a \in \R^\ell$,
$$
\|\sum_{i=1}^\ell a_i \xi_i\|_{\psi_2} \leq c\left(\sum_{i=1}^\ell a_i^2 \|\xi_i\|_{\psi_2}^2 \right)^{1/2}.
$$
\end{Lemma}

Clearly, a centred random vector $Y$ in $\R^\ell$ is $L$-subgaussian  if, for every $t \in \R^\ell$, $\|\inr{Y,t}\|_{\psi_2} \lesssim L \|\inr{Y,t}\|_{L_2}$. Also, if $Y$ is isotropic and $Y_1,...,Y_m$ are independent copies of $Y$, then for every $t \in \R^\ell$ and every $a \in \R^m$,
\begin{equation} \label{eq:sum-indepdent-psi-2}
\|\sum_{i=1}^m a_i\inr{Y_i,t}\|_{\psi_2} \leq c\left(\sum_{i=1}^m a_i^2 \|\inr{Y_i,t}\|_{\psi_2}^2 \right)^{1/2} \leq c^\prime L \|a\|_2 \|t\|_2.
\end{equation}

The proof that ${\cal A}$ is a large event is based on two facts. First, invoking a standard concentration argument (see, e.g. \cite{MR3837109}), it is straightforward to verify that for $m \gtrsim d$, with probability at least $1-2\exp(-c_0m)$,
$$
\lambda_{\max}(\Gamma_2/\sqrt{m}) \leq 1 + c^\prime(L)\sqrt{d/m}.
$$
Thus, if $m \geq c(L)d$ one can set $\kappa_1=2$.

\vskip0.3cm
The other component, concerning Equation \eqref{eq:rearrangement}, is treated in the next lemma.
\begin{Lemma} \label{lemma:A-large-uniform}
For every $L \geq 1$ there is a constant $C$ that depends on $L$ such that the following holds.
Let $X_1,...,X_m$ be independent copies of an isotropic, $L$-subgaussian random vector $X$ in $\R^d$. Then for $1 \leq k \leq m$, with probability at least $1-2\exp(-t^2)$,
$$
\sup_{\{a \in \R^m : \|a\|_2 \leq 1, \|a\|_0 \leq k\}}  \left\|\sum_{i=1}^m a_i X_i\right\|_2 \leq C(L) \left(\sqrt{d+k\log(em/k)} + t\right).
$$
\end{Lemma}

Consider $k=\theta m $ such that $k \log(em/k) \geq d$. Setting $t=\sqrt{k\log(em/k)}$ it follows that
$$
\sup_{\{a \in \R^m : \|a\|_2 \leq 1, \|a\|_0 \leq \theta m \}}  \left\|\sum_{i=1}^m a_i X_i\right\|_2 \leq c_2(L) \sqrt{\theta \log(e/\theta)} \sqrt{m}
$$
with probability at least $1-2\exp(-\theta m \log(e/\theta))$. Hence, in the context of \eqref{eq:rearrangement}, for any $0<\delta \leq 1/4$ one may select any $\theta \leq \delta^2/c_2^2\log(e/\delta)$.

Returning to Theorem \ref{thm:bernoulli-good}, note that $L$ is an absolute constant and that $\kappa_1=2$. Moreover, $\delta$ and $\theta$ are constants that depend on $L$, and $m=cn$, $d=c^\prime d^*(K)$ for constants $c$ and $c^\prime$ that also depend on $L$---which means they are all absolute constants. As a result, $Pr({\cal A}) \geq 1-2\exp(-c^{\prime \prime} d^*(K))$ and $\Gamma$ is an optimal Dvoretzky-Milman ensemble with probability at least $1-2\exp(-\tilde{c}d^*(K))$, as claimed.
\endproof
\vskip0.3cm
\noindent{\bf Proof of Lemma \ref{lemma:A-large-uniform}.} Set $U_k=\{a \in \R^m : \|a\|_2 \leq 1, \|a\|_0 \leq k\}$ and note that
$$
\sup_{a \in U_k}  \left\|\sum_{i=1}^m a_i X_i\right\|_2 = \sup_{v \in B_2^d} \sup_{a \in U_k} \sum_{i=1}^m a_i \inr{X_i,v}.
$$
By convexity, it suffices to control the supremum over two $1/10$-nets: $\tilde{V} \subset B_2^d$ and  $\tilde{U} \subset U_k$. Given $I \subset \{1,...,m\}$, let $B_2^I$ be the set of points in $B_2^m$ that are supported on $I$ and observe that $U_k = \bigcup_{|I|=k} B_2^I$. Therefore, by a volumetric estimate, $|\tilde{V}| \leq \exp(c_0d)$ and $|\tilde{U}| \leq \exp(c_0k \log(em/k))$ for an absolute constant $c_0$. Fix $v \in \tilde{V}$ and $u \in \tilde{U}$ and by \eqref{eq:sum-indepdent-psi-2}, the random variable $y=\sum_{i=1}^m u_i \inr{X_i,v}$ satisfies that
$\|y\|_{\psi_2} \leq  cL \|v\|_2 \|u\|_2 \leq cL$.
Using the equivalent formulation of the $\psi_2$ norm from Lemma \ref{thm:psi}, we have that $Pr(|y| \geq c^\prime Lx) \leq 2\exp(-x^2)$. The proof is completed by setting $x=c_1\max\{\sqrt{d}, \sqrt{k\log(em/k)}\} + t$ followed by the union bound.
\endproof

\subsection{Proof of Theorem \ref{thm:log-concave}.}
Let $X$ be an isotropic, log-concave random vector in $\R^d$. By Borell's Lemma (see, e.g., \cite{MR3185453}), $X$ satisfies $\psi_1-L_2$ norm equivalence with an absolute constant $C$. In particular, for any $q>2$, $X$ satisfies an $L_q-L_2$ norm equivalence with constant $Cq$. Set, for example, $q=6$ and the required norm equivalence in Assumption \ref{ass:on-X} is satisfied for an absolute constant.

The assumption that all the one-dimensional marginals $\inr{X,v}$ and $\inr{Z,t}$ have densities that are uniformly bounded by $L^\prime$ implies that Assumption \ref{ass:small-ball-assumption} holds as well. Following the notation used in the two assumptions we denote the constants by $L$ and $\kappa_0$ respectively (and clearly $\kappa_0$ depends only on $L^\prime$).

As in the proof of Theorem \ref{thm:bernoulli-good}, the difficulty is in verifying that ${\cal A}$ has large probability. Thankfully, that has been addressed as part of the study of the extremal singular values of random matrices with iid log-concave rows (columns). The wanted estimate was established in \cite{MR2601042}, with a minor but important caveat from our perspective: an upper bound on $m$ as a function of $d$. That caveat was removed in \cite{MR3184689}, (see Thm~11.4.1 and Prop.~11.4.4):

\begin{Theorem} \label{Thm:A-k}
There exist absolute constants $c$ and $C$ such that the following holds. Let $Y$ be an isotropic, log concave random vector in $\R^\ell$,and let $Y_1,....,Y_m$ be independent copies of $Y$. Then for $1 \leq k \leq m$ and $u \geq 1$, with probability at least $1-c\exp(-u)$,
$$
\sup_{\{a \in \R^m : \|a\|_2 \leq 1, \ \|a\|_0 \leq k\}} \left\|\sum_{i=1}^m a_i Y_i\right\|_2 \leq C\left(u + \sqrt{k}\log(em/k) + \max_{1 \leq i \leq m} \|Y_i\|_2\right).
$$
Moreover, if $T$ is the random matrix whose rows are $(Y_i/\sqrt{m})_{i=1}^m$, then conditioned on the event $\max_{1 \leq i \leq m} \|Y_i\|_2 \leq (m\ell)^{1/4}$, with probability at least $1-c\exp(-\ell)-c\exp(-(m\ell)^{1/4})$, we have
$$
1-c^\prime \sqrt{\frac{\ell}{m}} \leq \lambda_{\min}(T) \leq \lambda_{\max}(T) \leq 1+c^\prime \sqrt{\frac{\ell}{m}}.
$$

\end{Theorem}

The second result we require here is due to Paouris \cite{MR2276533}. The formulation is from \cite{MR3185453}:
\begin{Theorem} \label{thm:Paouris}
There exists an absolute constant $C_1$ such that, for any integer $\ell$, any isotropic log-concave random vector $Y$ on $\R^\ell$ and any $u \geq 1$,
$$
Pr( \|Y\|_2 \geq C_1u\sqrt{\ell}) \leq 2\exp(-u\sqrt{\ell}).
$$
\end{Theorem}
The combination of Theorem \ref{Thm:A-k} and Theorem \ref{thm:Paouris} clearly implies the following:
\begin{Corollary} \label{cor:psi-phi-log-concave}
There exist absolute constants $c$ and $C$ such that, with probability at least $1-2\exp(-c \min\{(\ell m)^{1/4},\ell\})$, for every $1 \leq k \leq m$,
$$
\sup_{\{a \in \R^m : \|a\|_2 \leq 1, \ \|a\|_0 \leq k\}} \left\|\sum_{i=1}^m a_i Y_i\right\|_2 \leq C\left( (\ell m)^{1/4} + \sqrt{k}\log(em/k)\right),
$$
and
$$
1-c^\prime \sqrt{\frac{\ell}{m}} \leq \lambda_{\min}(T) \leq \lambda_{\max}(T) \leq 1+c^\prime \sqrt{\frac{\ell}{m}}.
$$

\end{Corollary}

\proof Fix $1 \leq k \leq m$ and set $u=(m/\ell)^{1/4} \geq 1$. By Theorem \ref{thm:Paouris}, with probability at least $1-2m\exp(-(m\ell)^{1/4}) \geq 1-2\exp(-c(m\ell)^{1/4})$, $\max_{1 \leq i \leq m} \|Y_i\|_2 \leq (m\ell)^{1/4}$. Therefore, applying Theorem \ref{Thm:A-k}, for any $1 \leq k \leq m$, it follows that with probability at least $1-2\exp(-c_1(m\ell)^{1/4})$,
\begin{equation} \label{eq:est-sparse-log-concave}
\sup_{\{a \in \R^m : \|a\|_2 \leq 1, \ \|a\|_0 \leq k\}} \left\|\sum_{i=1}^m a_i Y_i\right\|_2 \leq C\left(\sqrt{k}\log(em/k) + (m\ell)^{1/4}\right).
\end{equation}
The proof of the first claim is concluded by the union bound with respect to $1 \leq k \leq m$. The proof of the second claim is immediate.
\endproof

Finally, set $\ell=d$ and let $m \geq n$ to be specified in what follows. Consider the event from Corollary \ref{cor:psi-phi-log-concave}. By the second claim from that corollary we may set $\kappa_1=2$.

Now, given any $0<\delta<1/4$, let $m \geq (2C/\delta)^4 d^*(K)$ (and recall that one has the freedom to choose $m$ that is large enough). Set $\theta$ to satisfy that $C\theta \log (e/\theta) \leq \delta/2$ and therefore,
$$
\sup_{\{a \in \R^m : \|a\|_2 \leq 1, \ \|a\|_0 \leq \theta m\}} \left\|\sum_{i=1}^m a_i X_i\right\|_2 \leq \delta.
$$
Hence, it is evident that $Pr({\cal A}) \geq 1-2\exp(-c_2\min\{(nd)^{1/4},d\})$, as claimed.
\endproof

\subsection{A heavy-tailed example} \label{sec:example-heavy-tailed}
In the final example we present, the random vector $X$ can be heavy-tailed.
\begin{Theorem} \label{thm:general-example}
For $\kappa_0$, $L$ and $q>8$ there exist constants $c_0$ and $c_1$ that depend on $\kappa_0$, $L$ and $q$ such that the following holds.
Let  $K \subset \R^n$ be a convex body. Set $d = c_0d(K^*)$ and $m = c_1 n$. Assume that
\begin{description}
\item{$(1)$} $Z$ is an isotropic, $L$-subgaussian random vector in $\R^n$.
\item{$(2)$} $X$ is an isotropic random vector in $\R^d$ that satisfies $L_q$-$L_2$ norm equivalence with constant $L$ for some $q>8$.
\item{$(3)$} $\|X\|_2 \leq 100\sqrt{d}$ almost surely.
\item{$(4)$} $Z$ and $X$ satisfy the small-ball property with constant $\kappa_0$.
\end{description}

Then for $\beta \geq 1$, with probability at least $1-c_2m^{-\beta}$, we have that
$$
\sup_{v \in B_2^d} \|\Gamma_2\| \leq 2 \sqrt{m},
$$
and
$$
\sup_{\{\|a\|_2 \leq 1, \ \|a\|_0 \leq \theta m\}} \left\|\sum_{i=1}^m a_iX_i\right\|_2 \leq c_3\left(\sqrt{\beta} \sqrt{d} + \theta^{1/4}\sqrt{m}\right)
$$
for absolute constant $c_2$ and $c_3$.

In particular, with probability $1-c_2/m^\beta - 2\exp(-cd)$, $\Gamma$ is an optimal Dvoretzky-Milman ensemble.
\end{Theorem}

\begin{Remark}
Clearly, the constant $100$ in $(3)$ is arbitrary. Moreover, the assumption that $X$ has a bounded diameter almost surely can be relaxed considerably: all that is needed is that the event $\{ \max_{1 \leq i \leq m} \|X_i\|_2 \leq (md)^{1/4}\}$ has nontrivial probability. For the sake of simplicity we shall not pursue this further and keep the bounded diameter assumption.
\end{Remark}

The proof of Theorem \ref{thm:general-example} has appeared implicitly in \cite{MR3191978}, where it was shown that random matrix $\Gamma_2$ satisfies a Bai-Yin estimate: with probability at least $1-c/m^\beta$,
$$
\sup_{v \in S^{d-1}} \left|\frac{1}{m}\sum_{i=1}^m \inr{X_i,v}^2 -1 \right| \leq c^\prime \sqrt{\frac{d}{m}}.
$$
Clearly, on that event $\sup_{v \in B_2^d} \|\Gamma_2v\|_2 \leq 2\sqrt{m}$.

Moreover, the crucial part in the proof of the Bai-Yin estimate from \cite{MR3191978} is an upper estimate on
$$
\sup_{\{\|a\|_2 \leq 1, \ \|a\|_0 \leq k\}} \left\|\sum_{i=1}^m a_iX_i\right\|_2
$$
for a specific choice of $k$. A careful examination of the proof (see Lemma~2.3 and Theorem~2.4 in \cite{MR3191978}) reveals that the choice of $k$ is actually arbitrary:

\begin{Lemma} \label{lemma:general-est}
Following the notation of Theorem \ref{thm:general-example}, for $1 \leq k \leq m$ and $\beta \geq 1$, with probability at least $1-c/m^{\beta}$,
$$
\sup_{\{\|a\|_2 \leq 1, \ \|a\|_0 \leq k\}} \left\|\sum_{i=1}^m a_iX_i\right\|_2 \leq c^\prime\left(\sqrt{\beta}\max_{1 \leq i \leq m }\|X_i\|_2 + (mk)^{1/4} \right);
$$
here $c$ is an absolute constant and $c^\prime$ depends on $L$ and $q$.
\end{Lemma}

Since the proof of Lemma \ref{lemma:general-est} is almost identical to the combination of Lemma~2.3 and Theorem~2.4 from \cite{MR3191978}, we shall not reproduce it here. The combination of the two facts imply that $Pr({\cal A}) \geq 1-c/m^\beta$, and the assertion of Theorem \ref{thm:main-intro} holds by invoking Theorem \ref{thm:general-example}.
\endproof

\section{Generic Chaining and Bernoulli processes} \label{sec:chaining}
The proofs of Theorem \ref{thm:upper} and Theorem \ref{thm:lower} are based on Talagrand's \emph{generic chaining} method, which was introduced as a way of controlling the suprema of random processes. We refer to Talagrand's invaluable manuscript \cite{MR3184689} for an extensive survey on generic chaining and its applications. Here, we present the basic notions we require, and only in the setup that is needed in what follows.
\vskip0.3cm
Let $V \subset (\R^\ell,\| \ \|)$. A collection of subsets of $V$, $(V_s)_{s \geq 0}$, is an \emph{admissible sequence} of $V$ if $|V_0|=1$ and for every $s \geq 1$, $|V_s| \leq 2^{2^s}$. The idea is that the sets $V_s$ are increasingly fine approximations of $V$, and that the $\gamma_2$ functional, defined below, captures some ``aggregate complexity" of $V$.

To define the $\gamma_2$ functional, let $\pi_sv$ be the nearest point to $v$ in $V_s$ with respect to the norm $\| \ \|$, and set $\Delta_s v =\pi_{s+1}v-\pi_sv$.
\begin{Definition} \label{def:gamma-2}
Let
$$
\gamma_2(V,\| \ \|)= \inf \sup_{v \in V} \left(\sum_{s \geq 0} 2^{s/2}\|\Delta_s v\|+ \|\pi_0 v\|\right),
$$
where the infimum is taken with respect to all admissible sequences of $V$.
\end{Definition}
The $\gamma_2$ functional is a natural object when trying to obtain upper bounds on random processes that exhibit a subgaussian behaviour with respect to the norm $\| \ \|$. For example, let $Y$ be a symmetric random vector in $\R^\ell$, and consider $\sup_{v \in V} \inr{Y,v}$, the supremum of the random process $v \to \inr{Y,v}$. By ``subgaussian behaviour of the random process" we mean that for every $p \geq 2$ and every $v \in \R^\ell$, $\|\inr{Y,v}\|_{L_p} \leq L \sqrt{p} \|v\|$. Writing the telescopic sum
$$
\inr{Y,v} = \sum_{s \geq 0} \inr{Y,\Delta_s v} + \inr{Y,\pi_0v}
$$
it follows that for  ``$s$-links" in all possible chains, it suffices to control the random variables $\{\inr{Y,\Delta_s v} : v \in V\}$. The latter is a collection consisting of at most $2^{2^s} 2^{2^{s+1}} \leq 2^{2^{s+2}}$ random variables. The uniform control over the $s$-links is easily obtained because $Pr( |\inr{Y,\Delta_sv}| \geq cLu 2^{s/2} \|\Delta_s v\| ) \leq 2\exp(-u^2 2^s)$, and the uniform control over all chains is obtained by taking the union bound over $s$.

This remarkably simple argument reveals why the $\gamma_2$ functional is a natural upper bound on the supremum of a subgaussian process. Talagrand discovered that, in fact, the \emph{reverse inequality} was true for \emph{gaussian processes}. In particular, in the context of the standard gaussian random vector in $\R^\ell$, he showed that:
\begin{Theorem} \label{thm:MM}
There exist absolute constants $c$ and $C$ such that, for every integer $\ell$ and  every $V \subset \R^\ell$,
$$
c\gamma_2(V,\| \ \|_2) \leq \E \sup_{v \in V} \sum_{i=1}^\ell g_i v_i \leq C\gamma_2(V,\| \ \|_2).
$$
\end{Theorem}
The upper bound is (in some equivalent formulation) due to Fernique, though this version was established by Talagrand. The lower bound is Talagrand's celebrated \emph{Majorizing Measures Theorem}. The proof of both parts can be found in \cite{MR3184689}.
\vskip0.3cm

Theorem \ref{thm:MM} implies that behaviour of the supremum of a gaussian process is determined by a metric invariant relative to a single metric. Unfortunately, obtaining a \emph{lower bound} on the supremum of the corresponding Bernoulli process is considerably harder. In what follows we denote by $(\eps_i)_{i=1}^\ell$ independent, symmetric $\{-1,1\}$-valued random variables and the Bernoulli process indexed by $V \subset \R^\ell$ is
$$
v \to \sum_{i=1}^\ell \eps_i v_i.
$$

It should be stressed that $\E \sup_{v \in V} \sum_{i=1}^\ell \eps_i v_i$ and $\E \sup_{v \in V} \sum_{i=1}^\ell g_i v_i$ need not be equivalent, as the example of $V=B_1^\ell$ shows. At the same time, it is straightforward to verify that there are absolute constant $c$ and $c^\prime$ such that for any $V \subset \R^\ell$,
$$
\frac{c}{\sqrt{\log \ell}}  \E \sup_{v \in V} \sum_{i=1}^\ell g_i v_i \leq \E \sup_{v \in V} \sum_{i=1}^\ell \eps_i v_i \leq c^\prime \E \sup_{v \in V} \sum_{i=1}^\ell g_i v_i,
$$
which is the best that one can hope for.  

Despite the logarithmic gap, we show in what follows that for the indexing sets $V$ we are interested in, the expected supremum of the Bernoulli process is actually \emph{equivalent} to the expected supremum of the gaussian process. To that end we require three features of the Bernoulli vector. Firstly, by \cite{MR1244666}, we have that for $(a_i)_{i=1}^\ell \in \R^\ell$
$$
\| \sum_{i=1}^\ell \eps_i a_i \|_{L_p} \sim \sum_{i \leq p} a_i^* + \sqrt{p} \left(\sum_{i  > p} (a_i^*)^2 \right)^{1/2},
$$
where $(a_i^*)_{i=1}^\ell$ denotes the nonincreasing rearrangement of $(|a_i|)_{i=1}^\ell$.
\vskip0.3cm

Secondly, Bernoulli processes satisfy a version of Sudakov's minoration, a fact that was established by Talagrand (see, e.g., \cite{MR1102015}). Here we use an equivalent formulation, noticed independently by R.~Lata{\l}a and the author.

\begin{Theorem} \label{thm:sudakov-reformulation}
There exist absolute constants $c_1$ and $c_2$ for which the
following holds. Let $U \subset \R^\ell$. For every $\delta \leq \sup_{v \in U}
\|v\|_{2}$ and $2 \leq p \leq \ell$, if $|U|
\geq \exp(p)$ and the set $\{\sum_{i=1}^\ell \eps_i v_i : v \in U\}$
is $c_1\delta\sqrt{p}$-separated in $L_p$, then
$$
\E \sup_{v \in U} \sum_{i=1}^\ell \eps_i v_i \geq c_2\delta
\sqrt{p}.
$$
\end{Theorem}

Finally, the supremum of a Bernoulli process exhibits a gaussian-like concentration, a fact that is also due to Talagrand (see, e.g. \cite{MR1102015}):

\begin{Theorem} \label{thm:Bernoulli-conc}
There exists and absolute constant $c$ such that, for $x>0$,
$$
Pr\left( \left| \sup_{v \in V} \sum_{i=1}^\ell \eps_i v_i - \E \sup_{v \in V} \sum_{i=1}^\ell \eps_i v_i \right| > x \right) \leq 2\exp\left(-c \frac{x^2}{\sup_{v \in V} \|v\|_2^2} \right).
$$
\end{Theorem}

\section{Proof of Theorem \ref{thm:upper}} \label{sec:upper}
The proof of Theorem \ref{thm:upper} is based on a chaining argument for a Bernoulli process indexed by a random set.

Let $\Gamma=\Gamma_1^* \Gamma_2^*$ and set $U$ to be a $1/2$ net in $B_2^d$ whose cardinality is at most $\exp(10d)$. The goal is to show that with the wanted probability,
$$
\max_{v \in U} \sup_{t \in K^\circ} \inr{\Gamma u,t} \leq C \ell(K)
$$
for a constant $C$ that depends on $\kappa_1$, $L$ and $q$. To that end, for $v \in S^{d-1}$ and $t \in K^\circ$ set
$$
W_{v,t} = \inr{\Gamma_1^* \Gamma_2^* v ,t} = \frac{1}{\sqrt{m}}\sum_{i=1}^m \inr{X_i,v}\inr{Z_i,t}.
$$
Let $(\eps_i)_{i=1}^m$ be independent, symmetric $\{-1,1\}$-valued random variables that are also independent of $(X_i,Z_i)_{i=1}^m$. Set
$$
W^\circ_{v,t} = \frac{1}{\sqrt{m}}\sum_{i=1}^m \eps_i \inr{X_i,v}\inr{Z_i,t}
$$
and by the symmetry of $Z$ and $X$ and the independence of $(\eps_i)_{i=1}^m$, $(X_i)_{i=1}^m$ and $(Z_i)_{i=1}^m$ the process
$$
\left\{(v,t) \to W_{v,t}, \ v \in U, t \in K^\circ\right\}
$$
has the same distribution as the process
$$
\left\{(v,t) \to W_{v,t}^\circ, \ v \in U, t \in K^\circ\right\}.
$$

Let $s_0$ be the smallest integer such that $2^{s_0} \geq 20d$ and for an integer $s$, let $k_s$ be the smallest integer such that $\gamma 2^{s+3} \leq  k \log(em/k)$ for a constant $\gamma$ to be specified in what follows. Let $(T_s)_{s \geq 0}$ be an optimal admissible sequence of $K^\circ$ with respect to the $\ell_2$ norm. As always, denote by $\pi_s t$ the nearest point to $t$ in $T_s$, set $\Delta_st=\pi_{s+1}t-\pi_st$ and observe that
$$
\gamma_2(K^\circ,\| \ \|_2) =\sum_{s \geq 0} 2^{s/2} \|\Delta_s t\|_2 + \|\pi_0t\|_2.
$$
Clearly,
$$
\sum_{i=1}^m \eps_i \inr{X_i,v}\inr{Z_i,t} = \sum_{s \geq s_0} \sum_{i=1}^m \eps_i \inr{X_i,v}\inr{Z_i,\Delta_st} + \sum_{i=1}^m \eps_i \inr{X_i,v}\inr{Z_i,\pi_{s_0}t}.
$$
Conditioned on $(X_i,Z_i)_{i=1}^m$, for any $v \in \R^d$, $t \in \R^n$ and $I \subset \{1,...,m\}$,
$$
\left(\E |W_{v,t}^\circ|^p \right)^{1/p}= \left\|W_{v,t}^\circ \right\|_{L_p(\eps)} \leq \frac{1}{\sqrt{m}}\sum_{i \in I} |\inr{X_i,v}\inr{Z_i,t}| + c\sqrt{p} \left(\frac{1}{m}\sum_{i \in I^c} |\inr{X_i,v}\inr{Z_i,t}|^2 \right)^{1/2},
$$
 for an absolute constant $c$. In particular, let $I=I_{v,t}$ be the union of the sets of indices corresponding to the $k$ largest values of $(|\inr{X_i,v}|)_{i=1}^m$ and the $k$ largest values of $(|\inr{Z_i,t}|)_{i=1}^m$. Here and throughout we abuse notation and write $\inr{Z_i,t}^*$ (resp. $\inr{X_i,v}^*$) for the $i$-th largest element in the non-increasing rearrangement of $(|\inr{Z_j,t}|)_{j=1}^m$ (resp. $(|\inr{X_j,v}|)_{j=1}^m$). Therefore,
\begin{equation} \label{eq:basic-chaining-Bernoulli}
\left\|W_{v,t}^\circ \right\|_{L_p(\eps)} \leq \|(\inr{X_i,v})\|_2 \left(\frac{1}{m}\sum_{i \leq 2k} (\inr{Z_i,t}^*)^2 \right)^{1/2} + c\sqrt{p} \left(\frac{1}{m}\sum_{i \geq k} (\inr{X_i,v}^* \inr{Z_i,t}^*)^2 \right)^{1/2}.
\end{equation}
The key observation is the following:
\begin{Theorem} \label{thm:estimate}
There is an event ${\cal B}$ of probability at most $2\exp(-d)$ such that on ${\cal A} \backslash {\cal B}$ the following holds:
\begin{description}
\item{$(1)$} For every $v \in U$, $\|(\inr{X_i,v})\|_2 \leq \kappa_1 \sqrt{m}$.
\item{$(2)$} For every $t \in K^{\circ}$, $v \in U$ and every $s \geq s_0$,
    $$
     \left(\frac{1}{m}\sum_{i \geq k_s} (\inr{X_i,v}^* \inr{Z_i,\Delta_st}^*)^2 \right)^{1/2} \leq c_1(q,L)\|\Delta_st\|_2.
    $$
\item{$(3)$} For every $t \in K^{\circ}$ and every $s \geq s_0$,
    $$
    \left(\frac{1}{m}\sum_{i \leq 2k_s} \left(\inr{Z_i,\Delta_st}^*\right)^2\right)^{1/2} \leq c_2(L)\frac{2^{s/2}}{\sqrt{m}} \|\Delta_s t\|_2.
    $$
\item{$(4)$} For every $t \in K^{\circ}$ and $v \in U$,
    $$
    \left(\frac{1}{m}\sum_{i \leq 2k_{s_0}} \left(\inr{Z_i,\pi_{s_0}t }^*\right)^2\right)^{1/2} \leq c_2(L)\frac{2^{s_0/2}}{\sqrt{m}}\sup_{t \in K^\circ} \|t\|_2,
    $$
    and
    $$
    \left(\frac{1}{m}\sum_{i \geq k_{s_0}} (\inr{X_i,v}^* \inr{Z_i,\pi_{s_0}t}^*)^2 \right)^{1/2} \leq c_1(q,L) \sup_{t \in K^\circ} \|t\|_2.
    $$
\end{description}
\end{Theorem}

The proof of Theorem \ref{thm:estimate} is based on several standard facts that are outlined in the next lemma.

\begin{Lemma} \label{lemma:monotone}
Let $(\xi_i)_{i=1}^m$ be independent copies of a centred random variable $\xi$. Set $1 \leq k \leq m$ and $u \geq 2$.
\begin{description}
\item{(a)} If $\xi \in L_q$ and $1 \leq r<q$, then with probability at least $1-2u^{-qk}\exp(-c_0k\log(em/k))$,
$$
\xi_k^* \leq eu \|\xi\|_{L_q} \left(\frac{m}{k}\right)^{1/r},
$$
where $c_0 \sim (q/r -1)$. In particular, if $\|\xi\|_{L_q} \leq L \|\xi\|_{L_2}$ then on that event,
$$
\xi_k^* \leq eLu \|\xi\|_{L_2} \left(\frac{m}{k}\right)^{1/r}.
$$
\item{(b)} If $\xi$ is $L$-subgaussian then with probability at least $1-2\exp(-u^2 k\log(em/k))$,
$$
\xi_k^* \leq c_1 L u \|\xi\|_{L_2} \sqrt{\log(em/k)}
$$
and
$$
\left(\sum_{i \leq k} (\xi_i^*)^2 \right)^{1/2} \leq c_1 L u\|\xi\|_{L_2}\sqrt{k \log(em/k)}
$$
for an absolute constant $c_1$.
\end{description}
\end{Lemma}

\proof
The proof of the first claim is an immediate outcome of a binomial estimate. Indeed,
\begin{align*}
& Pr\left(\exists I, \ |I|=k \ : \forall i \in I \ |\xi_i| \geq eu \|\xi\|_{L_q} \left(\frac{m}{k}\right)^{1/r} \right) \leq \binom{m}{k} Pr^k \left(|\xi| \geq eu \|\xi\|_{L_q} \left(\frac{m}{k}\right)^{1/r} \right)
\\
\leq & \left(\frac{em}{k}\right)^k \cdot u^{-kq} \left(\frac{k}{em}\right)^{kq/r} = u^{-kq} \exp\left( -\left(\frac{q}{r}-1\right) k \log\left(\frac{em}{k}\right) \right).
\end{align*}

The first part of $(b)$ follows from the second part because
$$
\xi_k^* \leq \left(\frac{1}{k} \sum_{i \leq k} (\xi_i^*)^2 \right)^{1/2}.
$$
To prove the second part, let
$$
\|x\|_{[k]} = \max_{|I| = k } \left(\sum_{i \in I} x_i^2\right)^{1/2}
$$
which is a norm on $\R^m$, and set $B$ to be the unit ball of its dual norm.
Let $x=(\xi_i)_{i=1}^m$ and denote by $G$ the standard gaussian random vector in $\R^m$. Since $(\xi_i)_{i=1}^m$ is a $c_0L$-subgaussian random vector, it follows from gaussian dominance (e.g., via the chaining mechanism, see \cite{MR3184689}) that for every norm $\| \ \|$ on $\R^m$ and any $1 \leq p<\infty$,
$$
\left(\E \|x\|^p\right)^{1/p} \lesssim L \left(\E \|G\|^p\right)^{1/p}.
$$
Hence, by the strong-weak inequality for a gaussian random vector (which follows from the gaussian concentration theorem), there is an absolute constant $c_1$ such that
$$
\left(\E \|x\|_{[k]}^p \right)^{1/p} \lesssim L \left(\E \|G\|_{[k]}^p \right)^{1/p} \leq c_1 L \left( \E \|G\|_{[k]} + \sqrt{p} \sup_{z \in B} \|z\|_2 \right).
$$
It is straightforward to verify that $\E \|G\|_{[k]} \lesssim \sqrt{k \log(em/k)}$ and that $\sup_{z \in B} \|z\|_2 =1$. The claim follows by setting $p=u^2k \log(em/k)$ and invoking Chebyshev's inequality.
\endproof

\noindent{\bf Proof of Theorem \ref{thm:estimate}.} Property $(1)$ holds on the event ${\cal A}$; therefore, it is enough to verify that Properties $(2)-(4)$ hold with high probability.

Recall that for every $t \in \R^n$, $\inr{Z,t}$ is $L$-subgaussian and for every $v \in \R^d$, $\|\inr{X,v}\|_{L_q} \leq L\|v\|_2=L$. Observe that
$$
\left| \left\{ \Delta_s t : t \in T\right\} \right| \leq 2^{2^{s}} \cdot 2^{2^{s+1}} \leq 2^{2^{s+2}}, \ \ {\rm and} \ \ |U| \leq \exp(10d).
$$
Fix $s \geq s_0$ and consider the two parts of Lemma \ref{lemma:monotone} for $k \geq k_s$:
Part (a) for every $\xi=\inr{X,v}$ for $v \in V$ and $r=1+q/2$ and Part (b) for every $\xi=\inr{X,\Delta_s t}$, $t \in K^\circ$. The probability with which the two parts hold for a pair $(v,t)$ is $1-2\exp(-c_0k \log(em/k))$

Hence, the two parts hold uniformly for every $v \in U$ and $\Delta_s t$ with probability at least
$$
1-2^{2^{s+2}} \exp(10d) \cdot 2 \exp(-c_0k \log(em/k)).
$$
Taking the union bound for $k \geq k_s$, both parts hold uniformly with probability at least
\begin{align*}
1-2^{2^{s+2}} \exp(10d) \exp(-c_1k_s \log(em/k_s)) \geq & 1-2^{2^{s+2}} \exp(10d) \exp(-c_1 \cdot (\gamma/2) 2^{s+3}) 
\\
\geq & 1-\exp(2^{s+2}),
\end{align*}
for a suitable choice of the absolute constant $\gamma$.

On that event, and by the isotropicity of $Z$ and $X$, there is an absolute constant $c$ such that
\begin{align*}
& \left(\frac{1}{m}\sum_{i \geq k_s} (\inr{X_i,v}^* \inr{Z_i,\Delta_st }^*)^2\right)^{1/2}  \leq cuL\|v\|_2 \|\Delta_st\|_2 \cdot \left(\frac{1}{m} \sum_{j \geq k_s} \left(\frac{m}{j}\right)^{2/r} \log(em/j)\right)^{1/2}
\\
\leq & c_2(q) L \|\Delta_s t\|_2,
\end{align*}
as required. Property $(2)$ now follows from the union bound on $s \geq s_0$ and recalling that $2^{s_0} \geq 20d$.

Property $(3)$ follows in a similar fashion, thanks to the second component in Part (b) of Lemma \ref{lemma:monotone} and the union bound for $s \geq s_0$.

The proof of Property $(4)$ is identical to the proofs of Properties $(2)$ and $(3)$, by considering the set of the $L$-subgaussian random variables $\{ \inr{Z, \pi_{s_0}t} : t \in T\}$ whose cardinality is at most $2^{2^{s_0}}$. We omit the standard details.
\endproof

\noindent{\bf Proof of Theorem \ref{thm:upper}.} Let ${\cal B}$ be the event from Theorem \ref{thm:estimate}.  Conditioned on ${\cal A} \backslash {\cal B}$ and invoking \eqref{eq:basic-chaining-Bernoulli}, there is an absolute constant $c_0$ such that the following holds. Let $s \geq s_0$, $w \geq 2$, and set $p=w^22^{s+3}$. Hence, for every $v \in U$ and $t \in K^\circ$,
\begin{align*}
& \|W_{v,\Delta_s t}^\circ\|_{L_p(\eps)}
\\
\leq &  \|(\inr{X_i,v})\|_2 \left(\frac{1}{m}\sum_{i \leq 2k_s} (\inr{Z_i,\Delta_s t}^*)^2 \right)^{1/2} + c_0w2^{(s+3)/2} \left(\frac{1}{m}\sum_{i \geq k_s} (\inr{X_i,v}^* \inr{Z_i,\Delta_st}^*)^2 \right)^{1/2}
\\
\leq & c_1(q,L,\kappa_1) \left( \sqrt{m} \cdot \frac{2^{s/2}}{\sqrt{m}} \|\Delta_s t \|_2 + w2^{s/2} \|\Delta_s t \|_2 \right) \leq 2c_1 w 2^{s/2} \|\Delta_s t\|_2,
\end{align*}
By Chebyshev's inequality,
\begin{equation} \label{eq:single-in-proof-upper}
Pr_\eps \left( \left|W_{v,\Delta_s t}^\circ\right| \geq 2ec_1 w 2^{s/2} \|\Delta_s t\|_2\right) \leq \exp(-w^2 2^{s+3}).
\end{equation}
Recall that $|\{\Delta_s t: t\in T\}| \leq 2^{2^{s+2}}$, that $|U| \leq \exp(10d)$, and the choice of $w$. By the union bound,  \eqref{eq:single-in-proof-upper} holds uniformly for every $t \in K^\circ$ and $v \in U$ as long as $s \geq s_0$.  And, by the union bound for $s \geq s_0$, there is an absolute constant $c_2$ such that with probability at least $1-2\exp(-c_2w^2d)$ with respect to $(\eps_i)_{i=1}^m$,
$$
\left| \sum_{s \geq s_0} W_{v,\Delta_s t}^\circ \right| \leq c_3(q,L,\kappa_1) w \sum_{s \geq s_0} 2^{s/2} \|\Delta_s t\|_2 \leq c_3 w \ell(K),
$$
where we have used the fact that
$$
\sum_{s \geq s_0} 2^{s/2} \|\Delta_s t\|_2 \leq \gamma_2(K^\circ,\| \ \|_2) \lesssim \ell(K).
$$

Also, setting $p=w^22^{s_0}$, it follows that for every $t \in K^\circ$ and every $v \in U$,
\begin{align*}
& \|W_{v,\pi_{s_0} t}^\circ\|_{L_p(\eps)} 
\\
\leq & \|(\inr{X_i,v})\|_2 \left(\frac{1}{m}\sum_{i \leq 2k_{s_0}} (\inr{Z_i,t}^*)^2 \right)^{1/2} + c_0w2^{s_0/2} \left(\frac{1}{m}\sum_{i \geq k_{s_0}} (\inr{X_i,v}^* \inr{Z_i,t}^*)^2 \right)^{1/2}
\\
\leq & c_4(q,L,\kappa_1) w 2^{s_0/2} \|\pi_{s_0}t\|_2 \leq c_4 w \sqrt{d} \sup_{t \in K^\circ} \|t\|_2 \leq c_4 w \ell(K),
\end{align*}
because $d \leq d^*(K)$. Hence, by the union bound, with probability at least $1-2\exp(-c_5w^2d)$ with respect to $(\eps_i)_{i=1}^m$, for every $t \in K^\circ$ and $v \in U$,
$$
|W_{v,\pi_{s_0}t}^\circ| \leq c_6  w \ell(K).
$$
To summarize, conditioned on ${\cal A} \backslash {\cal B}$, with probability at least $1-2\exp(-cw^2d)$ with respect to $(\eps_i)_{i=1}^m$, for every $v \in U$ and $t \in K^\circ$,
$$
|W_{v,t}^\circ| \leq |W_{v,\pi_{s_0}t}^\circ|+|\sum_{s \geq s_0} W_{v,\Delta_st}^\circ| \leq c^\prime(\kappa_1,q,L) w \ell(K).
$$
Now the claim follows using a suitable choice of $w=c^{\prime \prime}(\kappa_1,q,L)$ and a Fubini argument.
\endproof

\section{Proof of Theorem \ref{thm:lower}} \label{sec:lower}
Recall that $V$ is a $\rho$ net in $S^{d-1}$ for some fixed $0<\rho \leq 1/4$.
The idea behind the lower bound is that, mainly thanks to the small-ball property, the set $\Gamma_1 K^\circ$ consists of vectors in a ``good position". As a result, conditioned on a large subset of ${\cal A}$ and for every $v \in V$, the supremum of the Bernoulli process $W_{v,t}^\circ$ dominates $\sim \ell(K)$ with sufficiently high probability. That high probability leads to a uniform estimate in $V$.

\subsection{Preliminary estimates} \label{sec:prelim}
 The starting point is the following lemma, which was first noticed by A. Pajor.
\begin{Lemma} \label{lemma:gamma-to-net}
There exists an absolute constant $\kappa_2$ for which the following holds. For any $T \subset \R^n$ there is a subset $T^\prime \subset T$ of cardinality at most $\exp(\kappa_2 n)$ such that
$$
\E \sup_{t \in T^\prime} \sum_{i=1}^n g_it_i \geq \frac{1}{2} \E \sup_{t \in T} \sum_{i=1}^n g_it_i.
$$
\end{Lemma}

\proof Let $T^\prime$ be a maximal $\eps$-separated subset of $T$ whose cardinality is $\exp(k)$. By Sudakov's inequality (see, e.g., \cite{MR1102015,MR1036275}), there is an absolute constant $c$ such that
$$
\eps \leq \frac{c}{\sqrt{k}} \E \sup_{t \in T} \sum_{i=1}^n g_i t_i.
$$
Since $T \subset T^\prime + \eps B_2^n$ it follows that
\begin{align*}
\E \sup_{t \in T} \sum_{i=1}^n g_it_i \leq & \E \sup_{t \in T^\prime} \sum_{i=1}^n g_it_i + \E \sup_{v \in \eps B_2^n} \sum_{i=1}^n g_i v_i
\\
\leq & \E \sup_{t \in T^\prime} \sum_{i=1}^n g_it_i + \eps \sqrt{n} \leq \E \sup_{t \in T^\prime} \sum_{i=1}^n g_it_i + c\sqrt{\frac{n}{k}} \E \sup_{t \in T} \sum_{i=1}^n g_i t_i
\\
\leq & \E \sup_{t \in T^\prime} \sum_{i=1}^n g_it_i + \frac{1}{2} \E \sup_{t \in T} \sum_{i=1}^n g_i t_i
\end{align*}
provided that $k \geq 4c^2n$.
\endproof
With Lemma \ref{lemma:gamma-to-net} in mind, let $T^\prime \subset K^\circ$ satisfy that
\begin{equation} \label{eq:T-prime}
|T^\prime| \leq \exp(\kappa_2 n) \ \ \ {\rm and} \ \ \ \E \sup_{t \in T^\prime} \sum_{i=1}^n g_i t_i \geq \frac{1}{2} \ell(K)
\end{equation}
for a suitable absolute constant $\kappa_2$.

Next, let us describe the properties of $\Gamma_2^* V$ and of $\Gamma_1 T^\prime$ that are required in the proof. Let $s_0$ be the smallest integer that satisfies $2^s \geq 20d$ and set $k_s$ to be the smallest integer that satisfies $2^{s+3} \leq k \log(em/k)$.
\begin{Definition}
For $u \geq 2$ let $\Omega_u$ be the event on which:
\begin{description}
\item{$(1)$} For every $t \in \R^m$,
$$
\frac{1}{m} \sum_{i=1}^m \inr{Z_i,t}^2 \leq 2 \|t\|_2^2.
$$
\item{$(2)$} Let $(T_s)_{s \geq 0}$ be an optimal admissible sequence of $T^\prime$. Then for every $s \geq s_0$ and every $t \in T^\prime$,
    $$
    \left(\sum_{i \leq 2k_s} (\inr{Z_i,\Delta_s t}^*)^2\right)^{1/2} \leq Lu2^{s/2} \|\Delta_st\|_2;
    $$
    $$
    \inr{Z_j,\Delta_s t}^* \leq Lu \|\Delta_s t\|_2 \sqrt{\log(em/j)} \ \ \ {\rm for \ every } \ j \geq k_s;
    $$
and
    $$
    \left(\sum_{i \leq 2k_{s_0}} (\inr{Z_i,\pi_{s_0} t}^*)^2\right)^{1/2} \leq Lu2^{s_0/2} \|\pi_{s_0}t\|_2;
    $$
    $$
    \inr{Z_j,\pi_{s_0} t}^* \leq Lu \|\pi_{s_0} t\|_2 \sqrt{\log(em/j)} \ \ \ {\rm for \ every } \ j \geq k_{s_0}.
    $$
\item{$(3)$} $\sup_{\{\|a\|_2 \leq 1, \ \|a\|_0 \leq \theta m\}} \left\|\sum_{i=1}^m a_i X_i \right\|_2 \leq \delta \sqrt{m}$.
\item{$(4)$} For every $v \in V$ and every $j \geq d$,  $(\inr{X_i,v})^*_{j} \leq L u\left(\frac{em}{j}\right)^{q/q+2}$.
\item{$(5)$} For every $t_1,t_2 \in T^\prime$ and $v \in V$,
\begin{equation} \label{eq:cor-SB-uniform}
\left| \left\{ i : |\inr{X_i,v}| \geq \kappa_0 \ \ {\rm and } \ \  |\inr{Z_i,t_1-t_2} | \geq \kappa_0 \|t_1-t_2\|_2\right\} \right| \geq 0.98m.
\end{equation}
\end{description}
\end{Definition}

We begin by showing that $\Omega_u$ is a nontrivial event.

\begin{Theorem} \label{thm:estimate-lower}
There exist constants $c_0$ and $c_1$ that depend on $\kappa_0, \kappa_1$, $q$ and $L$ such that the following holds. Let
\begin{equation} \label{eq:cond-3}
m \geq c_0 \max\left\{\frac{d^*(K)}{\rho}, n\right\},
\end{equation}
and set $u \geq 4 \sqrt{\log(5/\rho)}$. There exists an event ${\cal B}_u$, such that $\Omega_u \supset {\cal A} \backslash {\cal B}_u$ and
$$
Pr({\cal B}_u) \leq 2\exp(-c_1 d \log u ).
$$
\end{Theorem}

The proof of Theorem \ref{thm:estimate-lower} is similar to the proof of Theorem \ref{thm:estimate}. We only outline the standard argument.

\proof Part $(1)$ holds with probability at least $1-2\exp(-c(L)m)$. Indeed, it suffices to consider a net in $S^{n-1}$ whose cardinality is at most $\exp(10n)$. By the $\psi_1$ version of Bernstein's inequality (see, e.g. \cite{MR3113826}) we have that
$$
Pr\left( \left|\frac{1}{m} \sum_{i=1}^m \inr{Z_i,t}^2 - \|t\|_2^2\right| \geq \frac{\|t\|_2^2}{4} \right) \leq \exp(-c_0(L)m),
$$
and the wanted estimate follows from the union bound and the choice of $m$.
\vskip0.3cm
\noindent The proof of Part $(2)$ is identical to the proof of Theorem \ref{thm:estimate}. Following the latter, it is evident that the claim holds with probability at least $1-2\exp(-c_2u^2d)$ provided that $u \geq c_3\sqrt{\log(5/\rho)}$ for constants $c_2$ and $c_3$ that depend on $L$ and $q$.
\vskip0.3cm
\noindent Part $(3)$ holds on the event ${\cal A}$.
\vskip0.3cm
\noindent To prove Part $(4)$, let $r=1+q/2$ and recall that $\|\inr{X,v}\|_{L_q} \leq L\|v\|_2$. By Lemma \ref{lemma:monotone}, with probability at least $1-u^{-qj}\exp(-c_3(q)j\log(em/j))$
\begin{equation} \label{eq:theta-m-cord}
(\inr{X_j,v})^* \leq L u \left(\frac{em}{j}\right)^{1/r}.
\end{equation}
Since $|V| \leq \exp(d\log(5/\rho))$ and $m \geq c_4 d/\rho$, it follows from the union bound that  with probability at least
$1-2\exp(-c_5(q) d (\log u + \log(5/\rho)))$, \eqref{eq:theta-m-cord} holds for every $v \in V$ and every $j \geq d$. Indeed, note that
$$
2^{s_0} \log(em/2^{s_0}) \geq 20 d \log(em/10d) \geq 20d\log(c_4/10\rho) \geq 20d\log(5/\rho)
$$
provided that $c_4$ is a sufficiently large absolute constant.
\vskip0.3cm
\noindent Finally, Part $(5)$ is evident by the small-ball assumption on each $\inr{Z,t}$ and $\inr{X,v}$. Indeed, for every $t \in \R^n$, $Pr(|\inr{Z,t}| \leq \kappa_0 \|t\|_2) \leq 1/1000$. By a binomial estimate there is an absolute constant $c_5$ such that with probability at least $1-2\exp(-c_5 m)$, $|\{i: |\inr{Z_i,t}| \geq \kappa_0 \}| \geq 0.99m$. A similar argument used for $|\inr{X,v}|$ and the fact that $|T^\prime| \leq \exp(\kappa_2 m)$ and $|V| \leq \exp(d\log(5/\rho))$ completes the proof---by invoking the union bound and recalling the choice of $m$.
\endproof

\subsection{A lower bound on Bernoulli processes } \label{sec:Bernoulli}

A crucial component in the proof of the lower bound is an equivalence result between the expected supremum of a Bernoulli process and of a gaussian one indexed by the same set---under certain structural assumptions on the indexing set.

\begin{Theorem} \label{thm:L-p-equiv-Bernoulli-Gauss}
For every $0<\lambda<1$ there exists a constant $c_0=c_0(\lambda)$ for
which the following holds. Let $U \subset \R^m$ and assume that there is $\eta>0$ such that
for every $v,w \in U$ and every $1 \leq p \leq \lambda m$,
$$
\|\sum_{i=1}^m \eps_i(v-w)_i\|_{L_p} \geq \eta \|\sum_{i=1}^m
g_i(v-w)_i\|_{L_p}.
$$
Then
$$
\E \sup_{v \in U} \sum_{i=1}^m \eps_i
v_i  \geq c_0\eta \E \sup_{v \in U} \sum_{i=1}^m g_i v_i.
$$
\end{Theorem}
The proof of Theorem \ref{thm:L-p-equiv-Bernoulli-Gauss} is almost identical to the proof of the majorizing measures theorem. It is based on Talagrand's construction of an admissible sequence, by showing
that the functional $\phi(U)=\E \sup_{v \in U} \sum_{i=1}^m \eps_i v_i $ satisfies the following growth condition:
\begin{Definition} \label{def:growth-condition}
Let $(U,d)$ be a metric space. The functional $\phi$ satisfies the growth condition if there are $r \geq 4$ and $c_0$ such that for every integer $s$, every $k=2^{2^s}$ and every $a>0$ the following holds. If $v \in U$, $v_1,...,v_k \in B(v,ra)$ are $a$-separated and $H_i \subset U \cap B(v_i,a/r)$ then
$$
\phi\Bigl(\bigcup_{\ell \leq m} H_\ell \Bigr) \geq c_0a2^{s/2} + \min_{\ell \leq m} \phi(H_\ell).
$$
\end{Definition}
The growth condition was used by Talagrand to construct an admissible sequence for the metric space $(U,d)$.
\begin{Theorem} \label{thm:admissible} \cite{MR3184689}
If $(U,d)$ satisfies the growth condition in Definition \ref{def:growth-condition} with parameters $r$ and $c_0$ then
$$
\gamma_2(U,d) \leq c_1  (\phi(U) + {\rm diam}(U,d)),
$$
where $c_1=c_1(r,c_0)$.
\end{Theorem}

\noindent{\bf Proof of Theorem
\ref{thm:L-p-equiv-Bernoulli-Gauss}.} By the majorizing measures theorem, $\gamma_2(U,\| \ \|_2) \sim \E \sup_{v \in U} \sum_{i=1}^m g_i v_i$, and therefore, by Lemma \ref{lemma:gamma-to-net}, it suffices to consider $U^\prime \subset U$ of cardinality at most $\exp(\kappa_2 m)$.
We will use Theorem \ref{thm:admissible} for $U^\prime \subset \R^m$ and upper bound $\gamma_2(U^\prime,\| \ \|_2)$.

Let us show that $U^\prime$ satisfies the growth condition for our choice of functional $\phi(H)=\E \sup_{v \in H} \sum_{i=1}^m \eps_i v_i$.

Fix $r>4$ to be named later. Let $k \leq \exp(\kappa_2m)$,
assume that $\{v_\ell, \  \ 1 \leq \ell \leq k\}$ is $a$-separated in
$\ell_2^m$, set $H_\ell \subset B(v_\ell,\sigma)$ for $\sigma=a/r$ and put $H = \cup_{\ell \leq k}
H_\ell$.

For every $v \in U^\prime$, let $X_v=\sum_{i=1}^m \eps_i v_i$, $Y_\ell=(\sup_{v \in H_\ell}
X_v)-X_{v_\ell}$ and $W=\max_{\ell \leq k} |Y_\ell - \E Y_\ell|$.
Applying the concentration inequality for Bernoulli processes, it follows that for
every $u>0$,
$$
Pr(|Y_\ell - \E Y_\ell | \geq u) \leq 2\exp(-c_1u^2/\sigma^2)
$$
for a suitable absolute constant $c_1$. Therefore, $Pr(W \geq u) \leq
2k\exp(-c_1u^2/\sigma^2)$, and thus $\E W \leq c_2\sigma \sqrt{\log
k}$.

Note that
\begin{equation} \label{eq:in-proof-MM-Bernoulli}
\sup_{v \in H} X_v \geq \max_{\ell \leq k} X_{v_\ell} + \min_{\ell \leq k} \E Y_\ell -W,
\end{equation}
and that
$$
\E \sup_{v \in H} X_v = \phi\Bigl(\bigcup_{\ell \leq k} H_\ell\Bigr) \ \ \ {\rm and} \ \ \ \E Y_\ell = \E\sup_{v \in H_\ell} X_v = \phi(H_\ell).
$$
Therefore, integrating \eqref{eq:in-proof-MM-Bernoulli}, there is an absolute constant $c_3$ such that
$$
\phi\Bigl(\bigcup_{\ell \leq k} H_\ell\Bigr) \geq \E \max_{\ell \leq k} X_{v_\ell} + \min_{\ell \leq k}
\phi(H_\ell) - c_3 \sigma \sqrt{\log k},
$$
and the growth condition follows once we show that 
$$
\E \max_{\ell \leq k} X_{v_\ell} \geq 2c_3\cdot (a/r) \sqrt{\log k}
$$ 
for an appropriate choice of $r$.

For $1 \leq \ell \leq k$, let $u_\ell = v_\ell - v_1$, and clearly
$$
\E\sup_{\ell \leq k} X_{v_{\ell}} = \E \sup_{\ell \leq k} X_{u_{\ell}}
$$
because $\E X_{v_1}=0$. Also, $u_i-u_j=v_i-v_j$, and thus $\{u_\ell : 1 < \ell \leq k\}$ is $a$-separated in $\ell_2^m$.

Let $G$ be the standard gaussian vector in $\R^m$. By the assumed equivalence between $\|X_{s-t}\|_{L_p}$ and $\|\inr{G,s-t}\|_{L_p}$ it is evident that for every $i \not = j, \ i,j > 0$ and every $2 \leq p \leq \lambda m$,
\begin{equation} \label{eq:lemma3.12}
\|X_{u_i-u_j}\|_{L_p} \geq \eta \|\inr{G,u_i-u_j}\|_{L_p} \geq c_4\eta\sqrt{p}\|u_i-u_j\|_2 \geq c_4 \eta \sqrt{p} a.
\end{equation}
Applying Sudakov's minoration for Bernoulli processes (Theorem \ref{thm:sudakov-reformulation}) for $p \sim \log k $,  we have that
$$
\E \sup_{\ell \leq k} X_{u_\ell} \geq c_5 \eta a \sqrt{\log k}.
$$
This concludes the proof of the growth condition for $k \leq \exp(\lambda m)$ by making a suitable choice of $r$ that is large enough.

Finally, since $k \leq \exp(\kappa_2 m)$, it remains to consider $\lambda m \leq \log k \leq
\kappa_2m$, and show that a similar norm equivalence to \eqref{eq:lemma3.12} is true in this range. Indeed, observe that for $\lambda m \leq p \leq m$,

\begin{align*}
\|X_{u_i-u_j}\|_{L_p} \geq & \|X_{u_i-u_j}\|_{L_{\lambda m}} \geq \eta \|\inr{G,u_i-u_j}\|_{L_{\lambda m}}
\\
\geq & c_6 \eta \sqrt{\frac{\lambda m}{p}}\|\inr{G,u_i-u_j}\|_{L_p} \geq
c_6 \eta \sqrt{\lambda}\|\inr{G,u_i-u_j}\|_{L_p}.
\end{align*}
Hence, the same estimate as in
\eqref{eq:lemma3.12} holds with a constant that depends only on $\lambda$, concluding the proof of the growth condition.
\endproof

\subsection{Returning to the lower bound}
With all the ingredients set in place, let $T^\prime$ be as in \eqref{eq:T-prime}, fix $v \in V$ and let us obtain a high probability lower bound on
$$
\sup_{t \in T^\prime} \frac{1}{\sqrt{m}} \sum_{i=1}^m \eps_i \inr{X_i,v} \inr{Z_i,t},
$$
conditioned on the event $\Omega_u$. Note that to be of any use, the lower bound should hold uniformly for every $v \in V$, implying that the required individual probability estimate has to ``defeat" $|V| \leq \exp(d\log(5/\rho))$.

\vskip0.3cm

For $v \in V$, let $J_v$ be the set of indices of the $\theta m$ largest coordinates of $(|\inr{X_i,v}|)_{i=1}^m$. Consider the random variables
$$
U_{v,t}^\circ = \frac{1}{\sqrt{m}} \sum_{i \in J_v^c} \eps_i \inr{X_i,v} \inr{Z_i,t},
$$
and the first order of business is to obtain a high probability lower bound on $\sup_{t \in T^\prime} U_{v,t}^\circ$, conditioned on $\Omega_u$.

Fix $v \in V$, set
$$
A_v= \left\{ \left(\frac{1}{\sqrt{m}} \inr{X_i,v} \inr{Z_i,t} \IND_{J_v^c}(i) \right)_{i =1}^m : t \in T^\prime\right\} \subset \R^m
$$
and observe that on $\Omega_u$,
\begin{align*}
\sup_{a \in A_v} \|a\|_2 \leq & \max_{t \in T^\prime} \left(\frac{1}{m}\sum_{i \in J_v^c} \inr{X_i,v}^2 \inr{Z_i,t}^2 \right)^{1/2}
\leq (\inr{X_{\theta m},v})^* \max_{t \in T^\prime} \left(\frac{1}{m}\sum_{i=1}^m \inr{Z_i,t}^2 \right)^{1/2}
\\
\leq & 2(\inr{X_{\theta m},v})^* \sup_{t \in K^\circ} \|t\|_2 \leq 2L u^2 \left(\frac{e}{\theta}\right)^{1/r} \sup_{t \in K^\circ} \|t\|_2.
\end{align*}

By the concentration theorem for Bernoulli processes, we have that
\begin{align*}
Pr_\eps \left( \left| \sup_{a \in A_v} \sum_{i=1}^m \eps_i a_i - \E_\eps \sup_{a \in A_v} \sum_{i=1}^m \eps_i a_i \right| \geq x \right) \leq & 2\exp\left(-c_0\frac{x^2}{\sup_{a \in A_v} \|a\|_2^2}\right)
\\
\leq & 2\exp\left(-c_1(L)\left(\frac{\theta^{2/r}}{u^4}\right)\frac{x^2}{\sup_{t \in K^\circ} \|t\|_2^2}\right),
\end{align*}
and setting $x = \frac{1}{2} \E_\eps \sup_{a \in A_v} \sum_{i=1}^m \eps_i a_i$, it is evident that
$$
\sup_{a \in A_v} \sum_{i=1}^m \eps_i a_i \geq \frac{1}{2}\E_\eps \sup_{a \in A_v} \sum_{i=1}^m \eps_i a_i
$$
with probability at least
$$
1-2\exp\left(-c_3(L) \left(\frac{\theta^{2/r}}{u^4}\right) \left(\frac{\E_\eps \sup_{a \in A_v} \sum_{i=1}^m \eps_i a_i}{\sup_{t \in K^\circ} \|t\|_2}\right)^2 \right).
$$
Therefore, if one can show that for a suitable constant $c_4$,
\begin{equation} \label{eq:lower-suggested}
\E_\eps \sup_{a \in A_v} \sum_{i=1}^m \eps_i a_i \geq c_4 \ell(K)
\end{equation}
and if
\begin{equation} \label{eq:cond-4}
c_3 c_4^2 \frac{\theta^{2/r}}{u^4} d^*(K) \geq 4 d \log(5/\rho),
\end{equation}
then with probability at least
\begin{equation} \label{eq:lower-prob-1}
1-2\exp\left(-\frac{c_3}{2} c_4^2 \theta^{2/r} d^*(K)\right)
\end{equation}
with respect to $(\eps_i)_{i=1}^m$, for every $v \in V$
\begin{equation} \label{eq:lower-est-1}
\sup_{t \in T^\prime}   \frac{1}{\sqrt{m}} \sum_{i \in J_v^c} \eps_i \inr{X_i,v} \inr{Z_i,t}  \geq \frac{c_4}{2} \ell(K).
\end{equation}
It should be noted that, following the restriction on $u$ from Theorem \ref{thm:estimate-lower}, a suitable choice of $u$ will be proportional to $\log(5/\rho)$.

\vskip0.3cm

For the proof of \eqref{eq:lower-suggested} we show that for every $v \in V$, the Bernoulli process restricted to the coordinates $J_v^c$ is likely to be large. That is based on Property $(5)$ satisfied by $\Omega_u$ --- most of the coordinates of $\inr{X_i,v}$ and $\inr{Z_i,t}$ are nontrivial.

\begin{Theorem} \label{thm:main-lower}
There is a constant $\kappa_3$ that depends only on $\kappa_0$ such that, conditioned on the event $\Omega_u$, for every $v \in V$,
$$
\E_\eps \sup_{a \in A_v} \frac{1}{\sqrt{m}} \sum_{i=1}^m \eps_i a_i = \E_\eps \sup_{t \in T^\prime} \frac{1}{\sqrt{m}} \sum_{i \in J_v^c} \eps_i \inr{X_i,v} \inr{Z_i,t} \geq \kappa_3 \ell(K).
$$
\end{Theorem}

\proof On $\Omega_u$ we have that for every $v \in V$ and $t_1,t_2 \in T^\prime$, the cardinality of
$$
I_{v,t_1,t_2} = \left\{i : |\inr{X_i,v}| \geq \kappa_0 \ \ {\rm and} \ \  \inr{Z_i,t_1-t_2}| \geq \kappa_0 \|t_1-t_2\|_2 \right\}
$$
is at least $0.98m$. Recall that $\theta <1/4$, implying that $|J_v| = \theta m \leq m/4$ and that 
$$
|I_{v,t_1,t_2} \cap J_v^c| \geq \frac{m}{2}. 
$$
Set
$$
I_v = \left\{i : |\inr{X_i,v}| \geq \kappa_0 \right\} \supset I_{v,t_1,t_2}.
$$
 By the contraction inequality for Bernoulli processes (see, e.g. \cite{MR1102015}) used twice,
\begin{align*}
\E_\eps \sup_{t \in T^\prime} \left|\frac{1}{\sqrt{m}} \sum_{i \in J_v^c} \eps_i \inr{X_i,v} \inr{Z_i,t} \right| \geq & \E_\eps \sup_{t \in T^\prime} \left|\frac{1}{\sqrt{m}} \sum_{i \in J_v^c \cap I_v} \eps_i \inr{X_i,v} \inr{Z_i,t} \right|
\\
\geq & \kappa_0 \E_\eps \sup_{t \in T^\prime} \left|\frac{1}{\sqrt{m}} \sum_{i \in J_v^c \cap I_v} \eps_i \inr{Z_i,t} \right|.
\end{align*}

Denote by $(a_i^*)$ the nonincreasing rearrangement of $(|a_i|)_{i \in J_v^c \cap I_v}$ and fix $p \leq m/4$. Invoking the characterization of the $L_p$ norm of a linear form of the Bernoulli vector, we have that, for any $t_1,t_2 \in T^\prime$
\begin{align*}
& \left\|\frac{1}{\sqrt{m}} \sum_{i \in J_v^c \cap I_v} \eps_i  \inr{Z_i,t_1-t_2} \right\|_{L_p(\eps)}
\\
\gtrsim &  \sum_{i \leq p} \inr{Z_i,t_1-t_2}^* + \sqrt{p} \left(\sum_{i > p} \left( \inr{Z_i,t_1-t_2}^*\right)^2 \right)^{1/2}
\\
\geq & \sqrt{p} \left(\sum_{i > p} \left( \inr{Z_i,t_1-t_2}^*\right)^2 \right)^{1/2} \gtrsim \sqrt{p} \kappa_0 \|t_1-t_2\|_2 \geq c_1 \kappa_0 \|\inr{G,t_1-t_2}\|_{L_p}.
\end{align*}
Hence, by Theorem \ref{thm:L-p-equiv-Bernoulli-Gauss} applied to the set
$$
F_v=\left\{ m^{-1/2}(\inr{Z_i,t}), \ i \in J_v^c \cap I_v, \ \ t \in T^\prime \right\}
$$
there is an absolute constant $c_2$ such that
$$
\gamma_2(F_v,\| \ \|_2) \leq c_2 \kappa_0^{-1} \E_\eps \sup_{t \in T^\prime} \frac{1}{\sqrt{m}} \sum_{i \in J_v^c \cap I_v} \eps_i \inr{Z_i,t}.
$$
At the same time, and considering only coordinates in $J_v^c \cap I_v$, it is evident that for every $x=m^{-1/2}(\inr{Z_i,t_1})$ and $y=m^{-1/2}(\inr{Z_i,t_2})$, $\|x-y\|_2 \geq c \kappa_0 \|t_1-t_2\|_2$. Indeed, $|J_v^c \cap I_{v,t_1,t_2}| \geq m/2$ and on $\Omega_u$, for every $i \in J_v^c \cap I_{v,t_1,t_2}$, $|\inr{Z_i,t_1-t_2}| \geq \kappa_0 \|t_1-t_2\|_2$. Thus,
$$
\gamma_2(F_v,\| \ \|_2) \geq c \kappa_0 \gamma_2(T^\prime,\| \ \|_2),
$$
and by the majorizing measures theorem, Lemma \ref{lemma:gamma-to-net} and the definition of $T^\prime$
$$
\E_\eps \sup_{t \in T^\prime} \frac{1}{\sqrt{m}} \sum_{i \in J_v^c} \eps_i  \inr{Z_i,v} \geq c_1(\kappa_0) \E \sup_{t \in T^\prime} \sum_{i=1}^m g_i t_i \geq \frac{c_1}{2}\ell(K),
$$
as required.
\endproof

Finally, to complete the proof of a high probability lower bound on $\sup_{t \in T^\prime} W_{v,t}^\circ$, we show that the contribution of each set of indices $J_v$ to the supremum of the Bernoulli process is not that big.

\begin{Lemma} \label{lemma:cutoff-large}
Let $u = 4\sqrt{\log(5/\rho)}$ and set $(X_i,Z_i)_{i=1}^m \in \Omega_u$. Then with probability at least $1-\exp(-u^2 d/2)$ with respect to $(\eps_i)_{i=1}^m$, for every $v \in V$
$$
\sup_{t \in K^\circ} \left|\frac{1}{\sqrt{m}} \sum_{i \in J_v} \eps_i \inr{X_i,v} \inr{Z_i,t} \right| \leq c(L)u\left(\delta+u\theta^{(q-2)/2(q+2)}\sqrt{\log(e/\theta)}\right) \ell(K).
$$
\end{Lemma}

\proof Recall that $s_0$ is the smallest integer such that $2^s \geq 20d$ and $k_s$ is the smallest integer such that $2^{s+3} \leq k  \log(em/k)$. Let $(T_s)_{s \geq 0}$ be an optimal admissible sequence of $K^\circ$ and fix $v \in V$. Denote by $I$ the union of the set of indices of the $k_s$ largest coordinates of $(\inr{X_i,v})_{i \in J_v}$ and of $(\inr{Z_i,t})_{i \in J_v}$. By the definition of $\Omega_u$,
\begin{align*}
& \left\|\frac{1}{\sqrt{m}} \sum_{i \in J_v} \eps_i \inr{X_i,v}\inr{Z_i,\Delta_st}\right\|_{L_p(\eps)}
\\
\lesssim & \frac{1}{\sqrt{m}} \sum_{i \in I} |\inr{X_i,v}\inr{Z_i,\Delta_st}| + \sqrt{p} \left(\frac{1}{m}\sum_{i \in J_v \backslash I} \inr{X_i,v}^2\inr{Z_i,\Delta_st}^2 \right)^{1/2}
\\
\lesssim & \left(\sum_{i \in J_v} \inr{X_i,v}^2\right)^{1/2}\left(\frac{1}{m}\sum_{i \leq 2k_s} (\inr{Z_i,\Delta_s t}^*)^2 \right)^{1/2} + \sqrt{p} \left(\frac{1}{m}\sum_{k_s \leq i \leq \theta m} \left(\inr{X_i,v}^*\inr{Z_i,\Delta_st}^*\right)^2 \right)^{1/2}
\\
\leq & c(L) u \|\Delta_s t\|_2 \left(\delta \sqrt{m} \cdot \frac{1}{\sqrt{m}} \sqrt{k_s \log(em/k_s)}  + \sqrt{p} \left(\frac{1}{m}\sum_{i=k_s+1}^{\theta m} \left(\frac{m}{i}\right)^{2/r} \log(em/i)\right)^{1/2} \right) =(*),
\end{align*}
where $r=1+q/2$. Setting $p=u^22^{s+3}$, we have that
\begin{equation} \label{eq:in-proof-1}
(*) \leq c_1(L) u \left(\delta  + u \theta^{1/2-1/r}\sqrt{\log(e/\theta)} \right) 2^{s/2} \|\Delta_s t\|_2,
\end{equation}
and by Chebyshev's inequality, followed by the union bound over all  $\{\Delta_s t: t \in K^\circ\}$, it is evident that with probability at least $1-2^{2^{s+2}}\exp(-u^22^{s+3}) \geq 1-\exp(-u^22^{s+2})$ with respect to $(\eps_i)_{i=1}^m$, \eqref{eq:in-proof-1} holds for every $t \in K^\circ$. Next, by the union bound for $s \geq s_0$, \eqref{eq:in-proof-1} holds with probability at least $1-2\exp(-u^22^{s_0+1})$ for every $t \in K^\circ$ and every $s \geq s_0$. On that event,

$$
\left|\frac{1}{\sqrt{m}} \sum_{i \in J_v} \eps_i \inr{X_i,v}\inr{Z_i,t-\pi_{s_0}t}\right| \leq c_1(L)u\left(\delta  + u \theta^{1/2-1/r}\sqrt{\log(e/\theta)} \right) \ell(K).
$$

Finally, by an identical argument, only this time for the random variables
$$
\frac{1}{\sqrt{m}} \sum_{i \in J_v} \eps_i \inr{X_i,v}\inr{Z_i,\pi_{s_0}t},
$$
for $p=u^22^{s_0+3}$,
\begin{align*}
& \left\|\frac{1}{\sqrt{m}} \sum_{i \in J_v} \eps_i \inr{X_i,v}\inr{Z_i,\pi_{s_0}t}\right\|_{L_p(\eps)}
\\
\lesssim & \left(\sum_{i \in J_v} \inr{X_i,v}^2\right)\left(\frac{1}{m}\sum_{i \leq 2k_{s_0}} (\inr{Z_i,\pi_{s_0}t}^*)^2 \right)^{1/2} + \sqrt{p} \left(\frac{1}{m}\sum_{ k_{s_0} \leq i \leq \theta m} \left(\inr{X_i,v}^*\inr{Z_i,\pi_{s_0}t}^*\right)^2 \right)^{1/2}
\\
\leq & c_2(L) u \left(\delta  + u \theta^{1/2-1/r}\log^{1/2}(e/\theta) \right) 2^{s_0/2} \|\pi_{s_0} t\|_2,
\end{align*}
and again, $2^{s_0/2} \|\pi_{s_0} t\|_2 \leq \sqrt{d} \sup_{t \in K^\circ} \|t\|_2 \leq \ell(K)$.

Chebyshev's inequality, the union bound over $\pi_{s_0}t \in T_{s_0}$, followed by the union bound over $v \in V$ completes the proof, recalling that $2^{s_0} \geq 20d$ and that $u = 4\sqrt{\log(5/\rho)}$.
\endproof

\vskip0.3cm
The proof of Theorem \ref{thm:lower} is now clear: one combines \eqref{eq:lower-est-1} (a lower bound that holds uniformly w.r.t $v \in V$ on the Bernoulli process for the indices in each $J_v^c$) and Lemma \ref{lemma:cutoff-large} (an upper bound that holds uniformly w.r.t $v \in V$ on the Bernoulli process for the indices in each $J_v$).

\endproof

\bibliographystyle{plain}
\bibliography{Dvor}

\end{document}